\documentclass[final]{siamltex}

\usepackage{epsfig}
\usepackage{amsmath}
\usepackage{amssymb}
\usepackage{subfigure}

\newtheorem{algorithm}{Algorithm}

\def\remark{{\sc Remark:\hspace{0.5em}}}

\def\eps{\varepsilon}
\def\ord{{\cal O }}
\def\cond{\mbox{cond\/}}

\def\range{\operatorname{range}}
\def\Zst{{\mathbb Z}}
\def\Cst{{\mathbb C}}

\def\Rst{{\mathbb R}}
\def\Cm{{\Cst^{2M+1}}}
\def\Cn{{\Cst^{2N+1}}}

\def\Psp{{\boldsymbol P}}
\def\PM{{\Psp_{\! M}}}
\def\QM{{P_M}}
\def\QN{{P_N}}
\def\PN{{\Psp_{\! N}}}
\def\PMj{{\Psp_{\! M_j}}}
\def\PMM{{\Psp^2_{\! M}}}
\newcommand{\pmlsp}{p^{(M)}}
\newcommand{\popt}{p_{*}}
\newcommand{\pnlsp}{p^{(N)}}
\newcommand{\pnstop}{p^{(\Nstop)}}

\newcommand{\TM}{T_{M}}
\newcommand{\VM}{V_{M}}
\newcommand{\Vm}{V_{M}}

\newcommand{\Vnull}{V_{0}}
\newcommand{\Vn}{V_{N}}
\newcommand{\Nopt}{N_{*}}
\newcommand{\Nstop}{N_{0}}
\newcommand{\Vnopt}{V_{\Nopt}}
\newcommand{\Vnopta}{V_{\Nopt}^{\ast}}
\newcommand{\VN}{V_{N}}

\newcommand{\Vna}{V_{N}^{\ast}}
\newcommand{\Vma}{V_{M}^{\ast}}
\newcommand{\Vni}{V_{N}^{+}}

\newcommand{\BKM}{B_{M}}
\newcommand{\BKMM}{B_{M+1}}
\newcommand{\betak}{\beta_{\ell-1}}
\newcommand{\betakk}{\beta_{\ell}}
\newcommand{\betakkk}{\beta_{\ell+1}}
\newcommand{\alphak}{\alpha_{\ell-1}}
\newcommand{\alphakk}{\alpha_{\ell}}

\newcommand{\cm}{c^{(M)}}
\newcommand{\cn}{c^{(N)}}
\newcommand{\cnstop}{c^{(\Nstop)}}
\newcommand{\bm}{b^{(M)}}

\newcommand{\smp}{s^{\eps}}
\newcommand{\noise}{\nu}

\newcommand{\TMM}{T_{M+1}}
\newcommand{\bmm}{b^{(M+1)}}
\newcommand{\cmm}{c^{(M+1)}}
\newcommand{\Ek}{E_{\ell}}
\newcommand{\Tk}{T_{\ell}}
\newcommand{\tk}{t^{(\ell)}}
\newcommand{\bk}{b^{(\ell)}}
\newcommand{\ck}{c^{(\ell)}}
\newcommand{\vk}{v^{(\ell)}}
\newcommand{\yk}{y^{(\ell)}}
\newcommand{\zk}{z^{(\ell)}}
\newcommand{\Tkk}{T_{\ell+1}}
\newcommand{\bkk}{b^{(\ell+1)}}
\newcommand{\ckk}{c^{(\ell+1)}}
\newcommand{\tkk}{t^{(\ell+1)}}
\newcommand{\vkk}{v^{(\ell+1)}}
\newcommand{\ykk}{y^{(\ell+1)}}

\newcommand{\bkl}{b_{\frac{\ell+1}{2}}}

\newcommand{\vkl}{v_{\frac{\ell+1}{2}}}

\newcommand{\bkm}{b_{-\frac{\ell+1}{2}}}

\newcommand{\vkm}{v_{-\frac{\ell+1}{2}}}

\newcommand{\bkn}{b_{-\frac{\ell}{2}}}
\newcommand{\vkn}{v_{-\frac{\ell}{2}}}
\newcommand{\Tcb}{\TM \cm = \bm}
\newcommand{\Tcbm}{\TMM \cmm = \bmm}
\newcommand{\Tcbk}{\Tk \ck = \bk}
\newcommand{\Tcbkk}{\Tkk \ckk = \bkk}
\newcommand{\lmin}{\lambda_{min}}
\newcommand{\lmax}{\lambda_{max}}
\newcommand{\copt}{c_{*}}
\newcommand{\coo}{c^{(0)}}

\newcommand{\Pnopt}{{\Psp}_{\!\Nopt}}

\newcommand{\In}{I_{N}}
\newcommand{\Inopt}{I_{\Nopt}}

\newcommand{\Tm}{T_{M}}

\def\conj#1{{\overline#1}}
\def\adots{\mathinner{\mkern1mu\raise1pt\vbox{\kern7pt\hbox{.}}\mkern2mu
   \raise4pt\hbox{.}\mkern2mu\raise7pt\hbox{.}\mkern1mu}}

\begin{document}

\title{\bf A Levinson-Galerkin algorithm for regularized trigonometric 
approximation\thanks{The author has been supported by project S7001-MAT, 
Schr\"odinger fellowship J01388-MAT of the Austrian Science foundation 
FWF, and NSF DMS grant 9973373.}}

\author{Thomas Strohmer\thanks{Department of Mathematics, 
        University of California, Davis, CA 95616-8633, USA;
        strohmer@math.ucdavis.edu.}}
\date{}
\maketitle

%
%

\begin{abstract}
Trigonometric polynomials are widely used for the approximation of a smooth 
function from a set of nonuniformly spaced samples. If the samples are 
perturbed by noise, a good choice for the polynomial degree of the 
trigonometric approximation becomes an essential issue to avoid overfitting 
and underfitting of the data. Standard methods for trigonometric least
squares approximation assume that the degree for the approximating
polynomial is known a priori, which is usually not the case in practice.
We derive a multi-level algorithm that recursively adapts to the 
least squares solution of suitable degree. We analyze under which
conditions this multi-approach yields the optimal solution. The proposed 
algorithm computes the solution in at most $\ord(rM + M^2)$ operations ($M$
being the polynomial degree of the approximation and $r$ being the number
of samples) by solving 
a family of nested Toeplitz systems. It is shown how the presented method
can be extended to multivariate trigonometric approximation. We demonstrate
the performance of the algorithm by applying it in echocardiography to 
the recovery of the boundary of the Left Ventricle of the heart.
\end{abstract}


\noindent

\pagestyle{myheadings}
\thispagestyle{plain}
\markboth{THOMAS STROHMER}{REGULARIZED TRIGONOMETRIC APPROXIMATION}

\section{Introduction}
\label{intro}

The necessity of recovering a function from a finite set of nonuniformly 
spaced measurements arises in areas as diverse 
as digital signal processing, geophysics, spectroscopy or medical imaging. 
The measurements $\{s_j\}_{j=1}^r$ are often distorted by several kinds of 
error. Hence a complete reconstruction of the function from the perturbed 
data $\smp_j = s_j+\noise_j$ is not possible. 
Often the function to be reconstructed is smooth, in which case 
a trigonometric polynomial of relatively low degree (compared to 
the possibly huge number of samples) can provide a good approximation to 
the function. This trigonometric approximation may be found by 
solving the least squares problem
\begin{equation}
\underset{p \in \PM }{ \min } \sum_{j=1}^{r} |p(x_j) -\smp_j|^2 w_j,
\label{LSP}
\end{equation}
where $w_j >0$ are weights and $\PM$ is the space of 
trigonometric polynomials of degree less than or equal to $M$.

Many efficient algorithms have been developed to solve \eqref{LSP}, e.g., see 
the articles \cite{New70,Dem89,RAG91,FGS95,Fas97}.
But surprisingly little attention has been paid to the problem of how to 
control the smoothness of the approximation in order to avoid overfitting
and underfitting of the data. An adaptation of the smoothness of the
approximation can be achieved of instance by providing a
suitable upper bound for the degree $M$ of the space $\PM$ in \eqref{LSP}.
In most of the aforementioned algorithms a necessary requirement to
get useful results in applications is that a good a priori guess of the 
degree of the trigonometric approximation is available. 
However a priori it is not clear what is a suitable degree for the
solution, in terms of how to choose a reasonable degree $M$
when solving \eqref{LSP}. Determining $M$ by ``trial and
error'' is certainly not a satisfactory alternative.

\begin{figure}[h]
\begin{center}
\subfigure[Original signal and perturbed samples]{
\epsfig{file=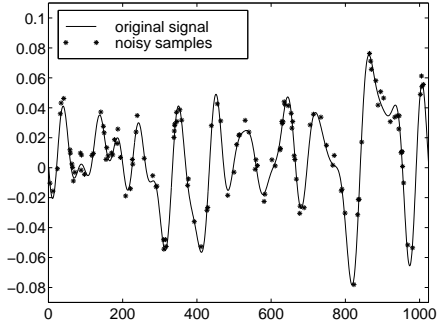,width=60mm}} \quad
\subfigure[Least squares approximation using a too large upper bound for
the degree of the polynomial (overfit)]{
\epsfig{file=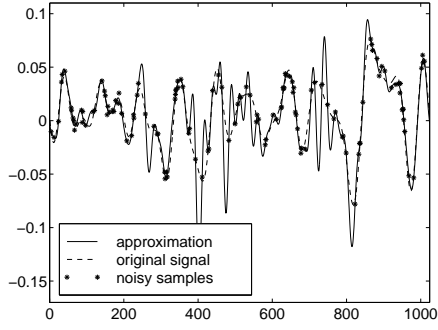,width=60mm}} \\
\subfigure[Least squares approximation using a too small upper bound for
the degree of the polynomial (underfit)]{
\epsfig{file=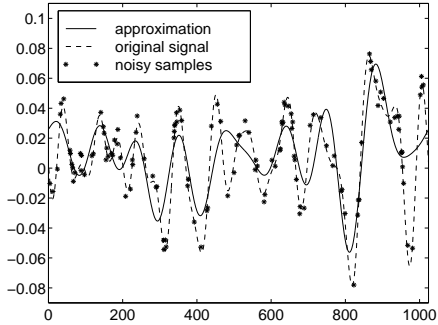,width=60mm}} \quad
\subfigure[Regularized approximation by proposed Levinson-Galerkin algorithm]{
\epsfig{file=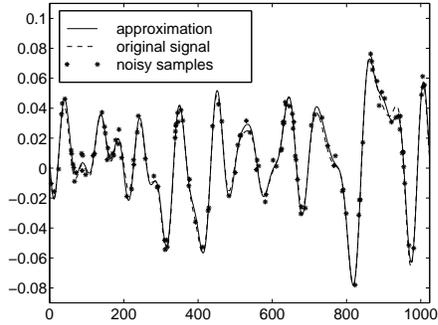,width=60mm}}
\caption{Controlling the smoothness of the solution is essential for 
trigonometric approximation from perturbed data in order to avoid
overfitting and underfitting of the data. The proposed Levinson-Galerkin 
algorithm automatically adapts to the least squares solution of optimal 
degree.}
\label{fig:degree}
\end{center}
\end{figure}

It is the goal of this paper to derive an efficient algorithm that 
computes the trigonometric approximation which provides the ``optimal'' balance 
between fitting the given data and preserving smoothness of the solution. 
Here optimality is meant in the sense that the solution has minimal 
degree among all trigonometric polynomials that satisfy a certain least 
squares criterion. The algorithm recursively adapts to the least squares
approximation of optimal degree by solving a family of nested Toeplitz 
systems in at most $\ord(Mr +M^2)$ operations,

If the data $\{\smp_j\}_{j=1}^{r}$ were
(i) unperturbed and (ii) stem from sampling a trigonometric polynomial
(with degree less than $r/2$), then the solution of \eqref{LSP} would
automatically have the appropriate degree, since the original function
could be completely recovered in this case. However the assumptions
(i) and (ii) are rarely met in applications and controlling the
smoothness of the solution becomes essential to avoid 
overfitting and underfitting of the data. 
If we choose the upper bound for the degree in \eqref{LSP} too large, 
the solution will almost always take on the maximal
possible degree, hence being too wiggly
and picking up too much noise (overfit), see also Figure~\ref{fig:degree}
(a)--(b). In the extreme case $2M+1 = r$ we will get an
interpolating polynomial, mostly with strong oscillations and far away from 
approximating the function between
the given samples. On the other hand, if we choose  $M$ too small, then
the approximation will be very smooth but poorly fitting the given data 
(underfit). Figure~\ref{fig:degree}(c) illustrates this behavior. 
The ``regularized'' trigonometric approximation obtained by the algorithm 
proposed in this paper --  to which we will refer as {\it Levinson-Galerkin 
algorithm} -- is shown in Figure~\ref{fig:degree}(d).

The paper is organized as follows. In Sections~\ref{s:multi} 
and~\ref{s:levgal} we
present the main results, including the Levinson-Galerkin algorithm
and a theoretical analysis that clarifies under which conditions this
algorithm provides optimal results. In Section~\ref{s:pred} we show
how properly chosen weights can be used as simple but efficient
tool to precondition the least squares problem.
Some aspects of extending the algorithm to multivariate
trigonometric polynomials are discussed in \ref{s:higherdim}.
In Section~\ref{s:curve} we present some applications in echocardiography.

Before we proceed we introduce some notation and conventions.
The inner product is denoted by $\langle \cdot , \cdot \rangle$, and the
conjugate transpose of a matrix $A$ by $A^{\ast}$.
The space of trigonometric polynomial of degree equal to or less than
$M$ is defined as
\begin{equation}
\PM = \left\{ p : p(x) = \sum_{k=-M}^{M} c_{k} e^{2\pi i kx } 
\right\}\,.
\label{pm}
\end{equation}
The norm of $p(x)=\sum_{k=-M}^{M} c_k e^{2\pi i kx} \in \PM$ is given by
\begin{equation}
\label{polynorm}
\|p\| = \left(\int \limits_{0}^{1} |p(x)|^2 \,dx\right)^{\frac{1}{2}}
= \left(\sum_{k=-M}^{M} |c_{k}|^2 \right)^{\frac{1}{2}}
= \|c\|\,,
\end{equation}
where $c=\{c_k\}_{k=-M}^M$.
In some applications it is advantageous to deal with complex-valued 
polynomials (see also Section~\ref{s:curve}), hence we do not restrict 
ourselves to the case of real-valued trigonometric approximation.

For $a=[a_{-M},a_{-M+1},\dots,a_{M-1},a_M] \in \Cm$ 
we define the orthogonal projections $\QN$ by
\begin{equation}
\label{defP}
\QN a = [0,\dots, 0, a_{-N}, a_{-N+1},\dots, a_{N-1}, a_N, 0, \dots,0]
\end{equation}
for $N=1,2,\dots,M$ and identify the image of $\QN$ with the 
$2N+1$-dimensional space $\Cst^{2N+1}$.

Let $p^{M}$ and $p^{N}$ be trigonometric polynomials of degree $M$ and
$N$ respectively, with coefficients vectors $c^M \in \Cst^{2M+1}, 
c^N \in \Cst^{2N+1}$.
If $N < M$, then we can always interpret $p^N$ as polynomial of degree
$M$ by adding to appropriate number of zero-coefficients 
and by doing so we are embedding the vector $c^N$ into a zero-padded vector of 
length $2M+1$. We will henceforth tacitly assume that such an
embedding has been made, when we compute expressions such as
$\|c^M - c^N\|$.

\section{Multi-level least squares approximation} \label{s:multi}


A standard method in numerical analysis to find the optimal balance between 
fitting the given data and preserving smoothness of the solution is to 
introduce a regularization parameter. The best value of this regularization 
parameter is then determined for instance by generalized cross 
validation~\cite{GHW79} or via the L-curve \cite{Han92}.
Here we understand regularization not as a way to stabilize ill-conditioned
problems, but in a broader context as a means of finding the best
compromise between fitting a given set of data and preserving smoothness
of the solution. As we will see, in our case it is not necessary to introduce 
an additional parameter,
since we can regularize the smoothness of the solution by varying the
parameter $M$ of the space $\PM$ in which we are searching for the solution 
of \eqref{LSP}. 

For the derivation of the algorithm we consider first the following situation.
Assume $\popt(x) =\sum_{k=-\Nopt}^{\Nopt} (\copt)_k e^{2 \pi i k x} \in \Pnopt$ and let 
$\smp=\{\smp_j\}_{j=1}^r, 2M+1 \le r$ with 
$\smp_j = s_j + \noise_j = \popt(x_j)+\noise_j$ be given noisy samples 
satisfying
\begin{equation}
\|\smp - s \|^2 \le \eps \|\smp\|^2.
\label{noisebound}
\end{equation}
For convenience we assume that $r$, the number 
of samples, is odd.

The aim is to approximate $\popt$ from the data $\{\smp_j\}_{j=1}^r$. 
Let us first assume that we already know that we are searching for
our least squares solution in the space $\Pnopt$.
In this case the coefficient vector of the polynomial that 
solves \eqref{LSP} is the least squares solution of
\begin{equation}
W \Vnopt c = W \smp \,,
\label{vanderwgt}
\end{equation}
where $\Vnopt$ is a $r \times (2\Nopt+1)$ Vandermonde matrix with entries
\begin{equation}
\label{vandermat}
(\VM)_{j,k}=e^{2\pi i k x_j}, \qquad j=1,\dots,r, \, k=-\Nopt,\dots,\Nopt
\end{equation}
and $W=\diag(\{\sqrt{w_j}\})$.

We will discuss the role and specific choice of the weights in more detail 
in Section~\ref{s:pred}. 
To reduce the notational burden we absorb the weight matrix $W$ in the 
Vandermonde matrix and in the sampling values. Thus for given degree, $M$
say, we consider the linear system of equations
\begin{equation}
\VM c = \smp \,
\label{vander}
\end{equation}
where $\VM$ is now the $r \times (2M+1)$ ``weighted Vandermonde'' matrix.
We will denote the least squares solution of \eqref{vander}
by $\cm=\{\cm_k\}_{k=-M}^{M}$ and the corresponding polynomial is
$\pmlsp(x) = \sum_{k=-M}^{M} \cm_k e^{2\pi ikx}$.

Since in general we do not know the optimal degree or {\em level} $M$ of the 
space $\PM$ in which we should solve the least squares problem, the situation 
becomes considerably more complicated. If we want to solve~\eqref{LSP} under 
the information \eqref{noisebound} without knowing the degree of the 
polynomial, one may argue that we have to accept any trigonometric polynomial
$p(x)=\sum_{k=-N}^{N} c_k e^{2\pi ikx}$ with
$\|\VN c - \smp \| \le \eps \|\smp\|$
as an approximate solution to $\popt$, since it is compatible with the only
knowledge we have on the data.

In general there may be infinitely many such polynomials, which raises the
questions of how to find a polynomial $p$ that yields a small
approximation error $\|\popt - p \|$ and at the same time can be computed 
efficiently.

\subsection{A multi-level algorithm and an efficient stopping criterion} \label{ss:stop}

The heuristic considerations above suggest the following approach.
\begin{algorithm}
\label{alg}
Set $N=0$ and solve $\Vnull \coo = \smp$. If $\coo$ satisfies the condition
$\|\Vnull \coo - \smp\| \le \eps \|\smp\|$, take $\coo$ as solution. Otherwise
set $N=N+1$ and solve 
\begin{equation}
\VN \cn = \smp \,,
\label{lsq}
\end{equation}
until $\cn$ satisfies for the first time the stopping criterion
\begin{equation}
\|\VN \cn - \smp\| \le \eps \|\smp\|.
\label{stop}
\end{equation}
at some level $N=\Nstop$. Set $\cnstop=\cn$. The approximation to $\popt$ is
then $\pnstop(x)=\sum_{k=-M}^{M} \cnstop_k e^{2\pi i kx}$.
\end{algorithm}

The stopping criterion~\eqref{stop} is well-defined, since it is definitely 
satisfied for $N=(r-1)/2$, in which case the left side in \eqref{stop} 
equals 0. Thus Algorithm~\ref{alg} selects among all least squares
solutions $\pnlsp, N=0,\dots,(r-1)/2$ that polynomial with minimal degree.

Algorithm~\ref{alg} and stopping criterion~\eqref{stop} 
can be justified by the following theoretical considerations.

One readily verifies that the matrices $\VN, N=0,\dots,(r-1)/2$ satisfy 
the relations:
\begin{itemize}
\setlength{\itemsep}{-0.5ex}
\setlength{\parsep}{-0.5ex}
\item[(i)] there exists a left-inverse $\Vni$ such that
\begin{equation}
\Vni \Vn = \In ,\quad \text{with} \,\, \Vni=(\Vna \Vn)^{-1} \Vna ,
\label{leftinv}
\end{equation}
where $\In$ is the identity matrix on $\Cn$.
\item[(ii)]
Let $a \in \Cm$ be the coefficient vector of some $p \in \PM$. Then
\begin{equation}
\Vn a = \Vm a \qquad \text{for all}\,\, N > M,  \, a \in \Cm.
\label{levels}
\end{equation}
\end{itemize}
In (ii) we have made use of the fact that the coefficient vector
$a$ can be interpreted as coefficient vector of 
a polynomial of degree $N$ by extending it to a vector of
length $2N+1$ via zero-padding. The matrix-vector multiplication
$\Vm a$ and equation~\eqref{levels} should be understood in this sense.


\begin{lemma}
\label{lemmastopactive}
If $N \ge \Nopt$ then $\cn$ satisfies $\|\Vn \cn - \smp \| \le \eps \|\smp\|$, 
hence stopping criterion~\eqref{stop} always becomes active at some level 
$\Nstop \le \Nopt$.
\end{lemma}
\begin{proof}
Note that $\Vn \Vni$ is the orthogonal projection into $\range (\Vn)$ and
$s \in \range (\Vnopt) \subseteq \range (\Vn)$ for $\Nopt \le N$, hence 
$\Vn \Vni s = s$. Therefore
\begin{align}
\|\Vn \cn - \smp \|^2 = & \|\Vn \Vni (s+\noise) - (s+\noise)\|^2
=\|\noise- \Vn \Vni \noise \|^2 \\
= & \|\noise \|^2 - \|\Vn \Vni \noise\|^2 \le \eps^2 \|\smp\|^2, 
\label{stopactive}
\end{align}
where we have used condition~\eqref{noisebound} in the last step.
It follows from \eqref{stopactive} that Algorithm~\ref{alg} 
terminates at some level $\Nstop \le \Nopt$.
\end{proof}

The following lemma shows that from the viewpoint of numerical stability
it is advisable to keep the level $N$ of the space $\PN$ in which we search for
our solution as small as possible.
\begin{lemma}
\label{lemmacond}
$\cond (\Vna \Vn) \ge \cond(\Vma \Vm)$ for $N\ge M$.
\end{lemma}
\begin{proof}
Since 
$$\QM (\Vna \Vn) \QM = \Vma \Vm \quad \text{for}\,\,M\le N, $$
Cauchy's Interlace Theorem~\cite{GL96} implies that $\cond (\Vna \Vn) \ge \cond(\Vma \Vm)$
for $N\ge M$.
\end{proof}

In the sequel we demonstrate that the fact that Algorithm~\ref{alg} terminates 
at some level $\le \Nopt$ is a desired property in many cases. 
We show that stopping criterion~\eqref{stop} is even optimum in a number of cases.

Let us first consider two special cases: (i) noisefree samples and (ii) uniformly 
spaced samples.

\subsubsection{Noisefree samples}
\label{noisefree}

Any reasonable stopping criterion has to satisfy the following lemma.
\begin{lemma}
\label{lemmanoisefree}
For noisefree data the stopping criterion~\eqref{stop} yields the exact solution.
\end{lemma}
\begin{proof}
One readily verifies that Algorithm~\ref{alg} terminates at level $\Nopt$.
Hence for $N =\Nopt$:
\begin{equation}
\|\pnlsp-\popt\| = \|\cn - \copt\| = \|\Vni \smp - \copt\| = 
\|\Vni \Vn \copt - \copt\| = 0,
\end{equation}
since $\Vni \Vn \copt = \copt$ for $N \ge \Nopt$.
\end{proof}

Lemmas~\ref{lemmanoisefree} and~\ref{lemmacond} together show that stopping 
criterion~\eqref{stop} yields the optimum solution for noisefree data
while providing maximum numerical stability.

{\subsubsection{Uniformly spaced samples}
\label{uniform}

If the sampling points $x_j, j=1,\dots,r$ are uniformly spaced and we choose
$w_j = 1/r$ as weights then a simple calculation shows that $\VN$ 
is unitary on $\Cn$, i.e., $\Vna \Vn =\In$ for $N =0,1,\dots,(r-1)/2$.

In this case 
\begin{equation}
\|\copt - \cn\| = \|\copt - \Vna \smp\|=\|\copt - \Vna \Vnopt \copt - \Vna \noise\|.
\label{unitary1}
\end{equation}
$N \ge \Nopt$ implies $\Vna \Vnopt =\In$ and hence
\begin{equation}
\|\copt - \cn\| = \|\Vna \noise\|.
\label{unitary2}
\end{equation}
Note that
\begin{equation}
\|\Vna \noise \| = \langle \Vna \noise, \Vna \noise \rangle = \langle \Vn \Vna \noise, \noise \rangle
= \|\Vn \Vna \noise \|,
\label{unitary3}
\end{equation}
since $\Vn \Vna$ is an orthogonal projection. Equation~\eqref{unitary3} yields 
\begin{equation}
\|\Vma \noise \| \le \|\Vna \noise \| \qquad \text{for} \,\,\, M \le N ,.
\label{unitary4}
\end{equation}
Consequently 
\begin{equation}
\|\copt -\cm \| \le \|\copt -\cn\| \qquad \text{if}\,\,\, \Nopt \le M \le N. 
\label{unitary5}
\end{equation}
Thus for uniformly spaced samples any stopping criterion should terminate Algorithm~\ref{alg}
at the latest at $N = \Nopt$. Under a mild condition on the coefficients $\copt$
we can show that the proposed stopping criterion provides the optimal solution among all least
squares solutions.

\begin{proposition}
\label{propo1}
Assume that the samples are regularly spaced. Then the solution $\pnstop$ computed via
Algorithm~\ref{alg} satisfies
\begin{equation}
\|\popt - \pnstop \| \le \|\popt -\pnlsp\| \qquad \text{for all}\,\,\,
 N \ge \Nopt.
\label{unitary6}
\end{equation}
If furthermore $\popt$ satisfies 
\begin{equation}
\label{decaycond}
\|(\Inopt-\QN) \copt \| \ge \|(\Inopt-\QN) \Vnopt^{\ast} \noise\|
\end{equation}
then
\begin{equation}
\|\popt - \pnstop \| \le \|\popt - \pnlsp \| \qquad \text{for all} \,\, N.
\label{unitary7}
\end{equation}
\end{proposition}
Condition~\eqref{decaycond} is satisfied e.g., if all coefficients of $\popt$ are 
larger than the relative noise level, i.e., $|c_k| \ge \eps \|\smp\|$.
\begin{proof}
Lemma~\ref{lemmastopactive} yields that $\Nstop \le \Nopt$, thus~\eqref{unitary5} 
implies~\eqref{unitary6}. 

To prove assertion~\eqref{unitary7} we only have to show that
$$\|\copt - \cnstop \| \le \|\copt - \cn\| \qquad \text{for all} \,\,N<\Nopt$$

For $N < \Nopt$ note that $(\copt-\Vna \Vnopt \copt)$ is orthogonal to 
$\Vna \noise$, since
\begin{align}
\langle \copt , \Vna \noise \rangle = & \langle \Vnopta \Vnopt \copt , \Vna \noise \rangle \\
= & \langle \Vnopt \copt , \Vnopta \Vna \noise \rangle 
= \langle \Vnopt \copt , \Vn \Vna \noise \rangle ,
\end{align}
hence
$$\langle \copt - \Vna \Vnopt \copt , \Vna \noise \rangle = 0.$$
Therefore
$$\|\copt - \cn \|^2 = \|\copt - \Vni \Vnopt \copt + \Vna \noise \|^2 
 = \|\copt - \Vni \Vnopt \copt\|^2 + \|\Vna \noise \|^2 \,.$$
In order to prove $\|\copt - \cnstop \| \le \|\copt - \cn\|$ for all $N < \Nopt$
we need to verify $ \|\copt - \Vni \Vnopt \copt\|^2 + \|\Vna \noise \|^2
\ge \|\Vnopta \noise\|^2$.
Since 
$$ \|\copt - \Vni \Vnopt \copt\|^2 = \sum_{|k|=N+1}^{\Nopt} |(\copt)_k|^2
= \|(\In - \QN) \copt\|^2$$
and
$$\|(\Vnopt)^{\ast} \noise \|^2 - \|\Vna \noise\|^2 = 
\sum_{|k|=N+1}^{\Nopt}|(\Vnopt^{\ast} \noise)_k|^2 = 
\|(\In - \QN) \Vnopta \noise\|^2,$$
the result follows now from the assumption~\eqref{decaycond}.
\end{proof}

\remark Proposition~\ref{propo1} shows that the least squares 
polynomial that gives the best approximation to $\popt$ is not
necessarily of degree $\Nopt$.

\subsubsection{Noisy nonuniform samples} \label{sss:general}

For noisy nonuniformly spaced data we observe that
$$\|\popt - \pnlsp\| = \|\copt - \cn\| \le 
\|\copt - \Vni \Vnopt \copt\| + \|\Vni \noise\|,$$
and for $N\ge \Nopt$
\begin{equation}
\label{general}
\|\popt - \pnlsp\| \le \|\Vni \noise\|,
\end{equation}
since $\|\copt - \Vni \Vnopt \copt\| =0$ in this case.

If $\VN$ is not unitary then $\|\Vni \noise\|$ is not necessarily monotonically 
increasing with increasing level $N$. One can argue heuristically that
since $\|\Vni\|$ is increasing with increasing level $N$ due to 
Lemma~\ref{lemmacond}, we may fairly assume that $\|\Vni \noise\|$ will
also increase (although not strictly monotonically). 
Also from the viewpoint of numerical stability it is reasonable to
keep the level $N$ small, since by Lemma~\ref{lemmacond} we
know that $\cond(\Vna \VN) \ge \cond(\Vma \VM)$ for $N \ge M$. This 
together with \eqref{general} suggests to choose a stopping criterion 
which terminates at or before level $\Nopt$, which is guaranteed for
stopping criterion~\eqref{stop} by Lemma~\ref{lemmastopactive}.

We can conclude that the stopping criterion
will provide excellent results if the noise level $\eps$ is small or if
the condition number of $\Vna \Vn$ is small (which implies that $\Vn$ is 
approximately unitary). In order to verify the latter it is useful
to have estimates for the condition number of $\Vna \Vn$. We will
address this issue in Proposition~\ref{propcond}.

\subsection{A Toeplitz system and trigonometric approximation} \label{ss:toep}

Instead of directly solving $\Vm \cm =\smp$ it is more efficient in our
case to consider the normal equations 
\begin{equation}
\Vma \Vma \cm = \Vma \smp .
\label{normal}
\end{equation}
The reason is that from a numerical point of view the structural properties of 
the matrix $\Vma \Vm$ are much more attractive than those of $\Vm$, which
in turn leads to faster numerical algorithms, see also Section~\ref{s:levgal}.

Set $\Tm = \Vma \Vm$ then a simple calculation shows that the entries of
the hermitian matrix $\Tm$ are 
\begin{equation}
(\Tm)_{k,l} = \sum_{j=1}^{r} w_j e^{2 \pi i (k-l) x_j},
\qquad k,l=-M,\dots,M.
\label{toeplitz}
\end{equation}
$\Tm$ is a Toeplitz matrix, since the entries $(\Tm)_{k,l}$ depend only on 
the difference $k-l$. Obviously $\Tm$ is invertible if $2M+1 \le r$. 

Following result is just a reformulation of \eqref{normal} together with
relation~\eqref{levels}, but since
it plays a key role in Section~\ref{s:levgal} it is helpful to state it 
in detail (cf.~also~\cite{Gro93a}).
\begin{theorem}
\label{thm1}
Given the sampling points $0 \le x_1 < \dots, x_r < 1$, samples 
$\{\smp_j\}_{j=1}^N$, positive weights $\{w_j\}_{j=1}^{r}$ and the degree 
$M$ with $2M+1 \le r$. The polynomial $\pmlsp \in \PM$ that solves 
\eqref{LSP} is given by
\begin{equation}
\pmlsp(x)=\sum_{k=-M}^M \cm_m e^{2 \pi i kx} \in \PM \,.
\label{lsppol}
\end{equation}
where its coefficients $\cm_k$ satisfy
\begin{equation}
\label{lspcoeff}
 \TM \cm = \bm \,\, \in \Cst^{(2M+1)^2},
\end{equation}
with 
\begin{equation}
\label{rightside}
\bm_k=\sum_{j=1}^{r} \smp_j w_j e^{2\pi i k x_j}
    \qquad \mbox{for $|k| \le M$}, 
\end{equation}
and $\Tm$ as defined in~\eqref{toeplitz}.
\end{theorem}

\section{Weights as simple preconditioner} \label{s:pred}

Vandermonde matrices are known to be ill-conditioned, if the nodes
$x_j$ are clustered~\cite{Gau91}. To improve the stability of the 
systems~\eqref{vander} and~\eqref{lspcoeff} we can use the weights as simple 
diagonal preconditioner. This leads to the problem of how to choose the
weights $w_j$.

We propose to use the size of the area of the Voronoi region~\cite{OBS92}
associated with the sampling point $x_j$ as weight $w_j$. In 1-D this 
reduces to
\begin{equation}
w_j = \frac{x_{j+1}-x_{j-1}}{2}\,.
\label{vorwgt}
\end{equation}
This choice is motivated by the following observations.

In this section we let $\VN$ denote the Vandermonde matrix defined
in~\eqref{vandermat} without weights.
Let $p \in \PM$ with coefficient vector $c$. Since
\begin{equation}
\langle \TM c,c \rangle = \langle W\VM c, W\VM c \rangle = \|W\VM c\|^2=
\sum_{j=1}^{r}|p(x_j)|^2 w_j, 
\label{smpenergy}
\end{equation}
the inequality
\begin{equation}
C_1 \|c\|^2 \le \sum_{j=1}^{r}|p(x_j)|^2 w_j \le C_2 \|c\|^2
\label{ineq}
\end{equation}
holds for all $c \in \Cst^{2M+1}$ with constants
$C_1=\lmin$ and $C_2=\lmax$, where $\lmin$ and $\lmax$ denote the 
minimal and maximal eigenvalue, respectively, of $\TM$. 

(i) The lower bound of $\TM$ is mainly determined by the large gaps  in
the sampling set. Suppose there is a large gap in the sampling set and
denote the interval corresponding to this gap by $\Gamma$
(hence $x_j \notin \Gamma$).
Choose a trigonometric polynomial $p\in \PM $ which, like the
prolate spheroidal functions, concentrates most of its energy in the
interval $\Gamma$.
Then the sampling values of $p$ will not pick up any information about the 
main concentration of the polynomial energy. Consequently if we use
no weights (or set $w_j=1$) we get
$$\sum_{j \notin \Gamma} |p(x_j)|^2 \ll \|p\|^2 =\|c\|^2.$$
For such a sampling set the lower frame
bound $C_1$ in the inequality~\eqref{ineq} must be small. 
Generically, large gaps and the ensuing lack of information always
results in bad condition numbers. This problem cannot be fixed by
preconditioning. 

(ii) On the other hand, we can choose a trigonometric polynomial that is
mainly concentrated in the region where the sampling points are located.
In this case the same local information is counted and added several times.
Thus 
$$\sum_{j \notin \Gamma}|p(x_j)|^2 \gg \|p\|^2=\|c\|^2$$ 
and the upper constant $C_2$ in \eqref{ineq} will be large. Yet, as
mentioned in (i) a cluster will not contribute much to the lower
bound and to the uniqueness of the problem. In this case the condition
number is large, because too much local information is given in
certain areas of the polynomial.

Problem (ii) can be addressed by introducing properly chosen weights.
The idea is to compensate for the local variation of the sampling
density by using weights in inequality~\eqref{ineq}.
Suppose that $0\leq x_1 \leq x_2 \leq \dots \leq x_r < 1$ is a
sampling set in $[0,1]$. Then a natural choice for the weights is $w_j =
(x_{j+1}-x_{j-1})/2$. Thus if many samples are clustered near a point
$x_j$, then the weight $w_j$ is small. If $x_j$ is the only sampling
point in a large neighborhood, then the corresponding weight is
large. This choice has not only been confirmed by extensive numerical
experiments~\cite{FGS95}, but also by the following optimization approach. 

A standard approach for the construction of preconditioners for a matrix 
$A$ is the following.  One attempts to find the matrix $P$ in a given class 
${\cal M}$ of matrices (e.g., the class of all circulant matrices or the 
class of all diagonal matrices) which solves 
\begin{equation}
\underset{P \in {\cal M}}{\min} \|I - P A \|_F \,,
\label{precond}
\end{equation}
where $\|.\|_F$ denotes the Frobenius norm.

In our setting this translates to the following optimization problem
\begin{equation}
\underset{W \in {\cal D}}{\min} \|I - (WV)^{\ast} WV \|_F \,,
\label{optim1}
\end{equation}
where ${\cal D}$ is the class of all $r \times r$ diagonal matrices 
and $I$ is the $(2M+1) \times (2M+1)$ identity matrix.

Note that we require that $w_j>0$ whereas~\eqref{optim1} could in principle 
yield weights that violate this condition. However since we will make
use of ~\eqref{optim1} in our actual algorithm, we are somewhat sloppy
here.

An alternative approach is to consider the solution of
\begin{gather}
\begin{array}{c}
\min \{\cond [(W V)^{\ast} WV]\} \\
\mbox{subject to $W \in {\cal D}$}.
\end{array}
\label{optim2}
\end{gather}
This optimization problem can be transformed to a general
eigenvalue problem, see~\cite{BEF94}, which can be solved by
convex optimization algorithms.

In the simple case of regular sampling it is easy to check that
the solution of both optimization problems is given
by $W=\diag(\{\sqrt{w_j}\})$ with $w_j = (x_{j+1}-x_{j-1})/2 = 1/r$. 
However in the more interesting case of nonuniform sampling neither 
problem~\eqref{optim1} nor \eqref{optim2} does in general have an 
analytic solution. Thus using these approaches for the actual
construction of a preconditioner would be ridiculous, since the
computational costs to solve these optimization problems are considerably
larger than solving the trigonometric approximation problem.
Nevertheless, solving~\eqref{optim1} and \eqref{optim2} numerically for
a variety of different examples is useful to get insight in the type of 
weights obtained by these approaches. 

The numerical results confirm the choice of the
Voronoi-type weights defined in~\eqref{vorwgt}. Sampling points in
densely sampled areas are assigned a small weight, whereas sampling points in
sparse sampled regions are assigned a large weight.
Two typical comparisons of the weights obtained via optimization and the
Voronoi weights are illustrated in Figure~\ref{fig:wgtopt}. 
In the first case we consider a sampling set with
high density at the endpoints and strongly decreasing density towards 
the center. The weights obtained by solving~\eqref{optim1} and
\eqref{optim2} are almost identical and are very close to the Voronoi
weights, as can be seen in Figure~\ref{fig:wgtopt}(a). The difference
at the endpoints is probably due to boundary effects.

In the second example we consider a random sampling set with 
several areas with high sampling density and relatively few samples
between these clusters. Again all three approaches give weights
that show a similar behavior, see Figure~\ref{fig:wgtopt}(b).
The condition number of the non-weighted Toeplitz matrix in this example
is 33, compared to the significantly smaller condition number 3.3 when
using Voronoi-type weights. Using the weights
obtained via \eqref{optim1} gives $\cond(\TM)=3.1$, and
for the weights resulting from \eqref{optim2} we get $\cond(\TM)=2.9$,
which is only a slight improvement compared to the Voronoi-type weights.

\begin{figure}
\begin{center}
\subfigure[]{
\epsfig{file=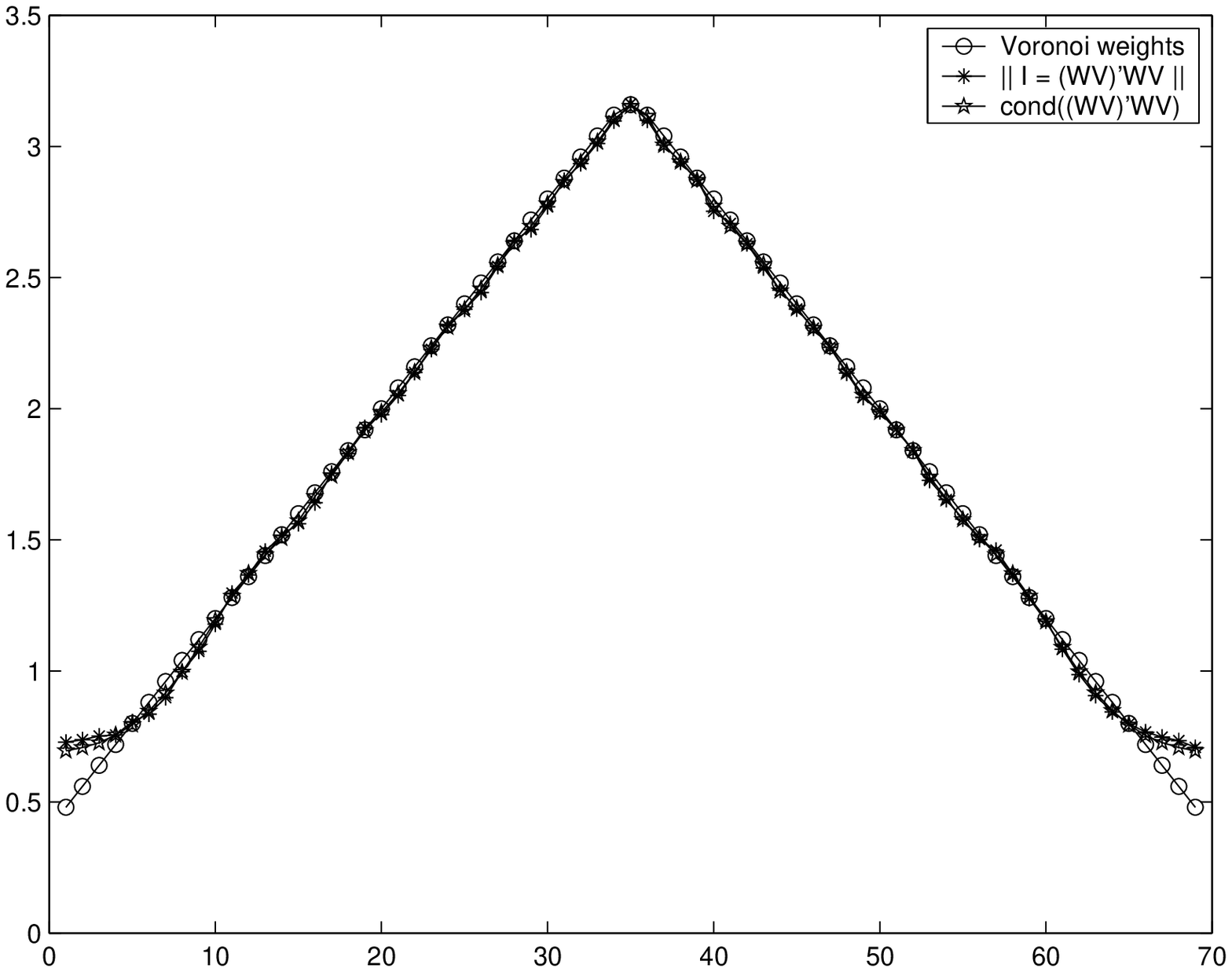,width=63mm}}
\subfigure[]{
\epsfig{file=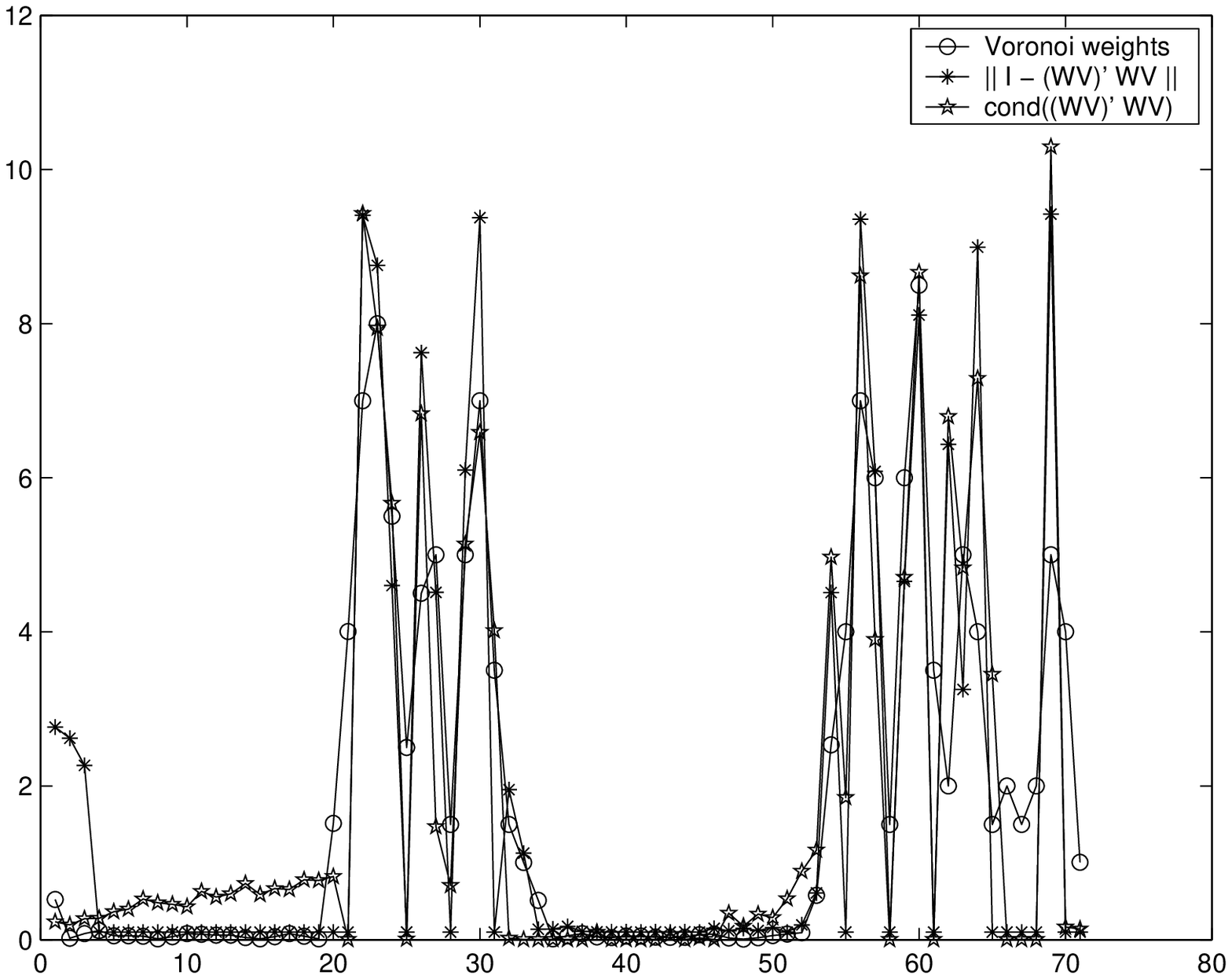,width=63mm}}
\label{fig:wgtopt}
\caption{Comparison of weights obtained by different approaches.}
\end{center}
\end{figure}

Obtaining good estimates for the condition number of a Toeplitz matrix is
a difficult problem. It is gratifying that by using the weights defined 
in~\eqref{vorwgt} it is possible to get an upper bound for the condition 
number.


\begin{proposition}[Gr\"ochenig, \cite{Gro93a}]
Assume that the sampling set $\{x_j\}_{j=1}^r$ satisfies
$$\max(x_{j+1}-x_j) := \gamma < \frac{1}{2M}$$ 
and set $w_j = (x_{j+1} - x_{j-1})/2$. Then the condition number of the 
Toeplitz matrix $\TM$ defined in~\eqref{toeplitz} is bounded by
\begin{equation}
\kappa(\TM) \le \left(\frac{1+\gamma}{1-\gamma} \right)^2 \,.
\end{equation}
\label{propcond}
\end{proposition}
 
\section{A Levinson-Galerkin algorithm for trigonometric approximation} 
\label{s:levgal} 

The method described in Algorithm~\ref{alg} can be seen as a Galerkin-type
approach, since we try to determine an approximation by searching for a 
solution in a finite-dimensional space spanned by orthogonal polynomials, 
and by increasing the dimension of the space we increase the resolution of 
our approximation by adding more and more details.

When we use Levinson's algorithm~\cite{GL96} to solve~\eqref{lspcoeff}
for $M=0,1,\dots,\Nstop$ the total computational effort would
be of $\ord(\Nstop^3)$, since the solution of each system $\Tm \cm =\bm$ 
requires $\ord(M^2)$ operations. Using one of the fast Toeplitz 
algorithms~\cite{AG88,CN96} reduces this effort to $\ord(k M \log M)$ 
for each level $M$, where $k$ is the number of iterations, thus leading
to a total of $\ord (k \Nstop^2 \log \Nstop)$ operations.
In this section we show that the systems $\Tm \cm =\bm, M=0,1,\dots,\Nstop$ 
can be solved in $\ord(\Nstop^2)$ operations and the total effort
(including the calculation of the entries of $\TM$ and the evaluation
of the stopping criterion~\eqref{stop}) for computing $\pnstop$ is
$\ord(r \Nstop +  \Nstop^2)$ operations.

Following observation is crucial for the derivation of the proposed 
Levinson-Galerkin algorithm.
\begin{lemma}\label{l:casc}
For fixed degree $M$ and $M+1$ let $\TM, \bm$ and $\TMM, \bmm$ be the 
Toeplitz matrices and right hand sides as defined in~\eqref{toeplitz} and
\eqref{rightside}, respectively. Then $\TM$ and $\bm$
are embedded in $\TMM$ and $\bmm$ in the following way:
\begin{gather}
\label{embed}
\TMM = 
\begin{bmatrix}
t_0 & \dots & \conj{t_{2(M+1)}} \\
\vdots & \fbox{$\TM$} & \vdots \\
t_{2(M+1)} & \dots & t_0
\end{bmatrix},\,\,
\bmm = 
\begin{bmatrix}
\setlength{\fboxsep}{4mm}
b_{-(M+1)}\\
\fbox{\!\!\!\!\!\!\!\!\!\! $\bm$\!\!\!\!\!\!\!\!\!}  \\
b_{M+1}
\end{bmatrix}. 
\end{gather}
\end{lemma}
\begin{proof}
\eqref{embed} follows immediately from the definition of $\TM$ and $\bm$ and
\eqref{levels}.
\end{proof}

Unfortunately the solutions $\cm$ and $\cmm$ of the systems
$\TM \cm = \bm$ and $\TMM \cmm = \bmm$ are not related is such a simple 
manner.  But we can exploit the nested structure of the family 
$\{\TM\}_{M=1}^{\Nstop}$ by solving the systems $\TM \cm = \bm$ recursively 
via a modified Levinson algorithm. The standard Levinson algorithm cannot
be applied directly, since it only addresses Toeplitz systems, where the 
principal leading sub-matrix and the principal leading sub-vector of the 
right hand side stay unchanged during the recursion, which is not the case
here. For $\TMM$ it does not matter, if we enlarge $\TM$ by appending
new entries below or above, whereas the right hand side $b^{(M)}$ cannot
be rearranged in such a way, the principal leading subvector of the
right hand side will be changed if we switch from $\bm$ at level $M$ to $\bmm$
at level $M+1$. 

To adapt Levinson's algorithm to our situation, we have to split
up the change from the system $\Tcb$ at level $M$ to the system $\Tcbm$ at
level $M+1$ into two separate steps. Instead of indexing the matrix
$\TM$ and the vectors $\bm,\cm$ by the degree $M$, it is therefore
advantageous to index them according to their dimension.
For clarity of presentation we reserve the subscript $(M)$ for the degree
of the polynomial and its coefficient vector respectively, and use the 
subscript $(\ell)$ when we refer to the 
dimension of the corresponding coefficient vector in $\Cst^{\ell}$. Thus
for even $\ell$, $\bk = [b_{-\frac{\ell}{2}+1}, \dots, b_{\frac{\ell}{2}}]^T 
\in \Cst^\ell$, and for odd $\ell$ we set
$\bk = [b_{-\frac{\ell-1}{2}}, \dots, b_{\frac{\ell-1}{2}}]^T$  
(whence $b^{(1)} = b_0$),
analogously for $\ck$. Further it is useful in the sequel to denote
$\tk = [t_1,\dots,t_\ell]^T$. Then the Toeplitz matrix $\Tk$ of size 
$\ell \times \ell$ is generated by the vector $[t_0, (t^{(\ell-1)})^T]^T$ with
$t_k = \sum_{j=1}^{r} w_j e^{2 \pi i k x_j}$ according to~\eqref{toeplitz}.

Assume we have already solved the system $\Tcb$ at level $M$ (with $\ell=2M+1$) 
and now we want to switch to the next level $M+1$. As we have agreed,
we do this in two steps. In the first step ($\ell \rightarrow \ell+1$) the
Toeplitz system can be written as
\begin{gather}
\begin{bmatrix}
\Tk         & \Ek \overline{t^{(\ell)}} \\
( \tk )^{T} \Ek  & t_0 
\end{bmatrix}
\begin{bmatrix}
\vk \\ 
\vkl
\end{bmatrix}
 =
\begin{bmatrix}
\bk \\
\bkl
\end{bmatrix}\,, 
\label{oddstep}
\end{gather}
where $\Ek$ is the rotated identity matrix on $\Cst^\ell$, i.e., 
$$
\Ek = 
\begin{bmatrix}
0    &         & 1 \\
     & \adots  &   \\
1    &         & 0 
\end{bmatrix} \,.
$$
System \eqref{oddstep} can be solved recursively by the standard Levinson 
algorithm \cite{Lev47,GL96}. To be more detailed, assume that we have already 
solved the system $\Tcbk$ for $\ell=2M+1$ and assume further that the 
solution of the $\ell$-th order {\em Yule-Walker system} $\Tk \yk = - \tk$ 
is available. Then the solution of~\eqref{oddstep} can be computed 
recursively by
\begin{align}
\vkl &=(\bkl - [\tk]^{T} \Ek \ck)/\betakk \notag \\
\vk &= \ck + \vkl \Ek \conj{\yk} \notag 
\end{align}
where
\begin{align}
\betakk &=t_0+[\tk]^{T}\conj{\yk}=(1-\alphak \conj{\alphak}) \betak 
\notag \\
\alphakk &= - (t_{\ell+1}+[\tk]^{T} \Ek \yk)/\betakk \notag \\
\zk &= \yk + \alphakk \Ek \conj{\yk} \notag  \\
\ykk &= \begin{bmatrix} \zk \\ \alphakk \end{bmatrix}\notag \,.
\end{align}

Now we can proceed to the second step ($\ell+1 \rightarrow \ell+2 = 2(M+1)+1$), 
where the Toeplitz system can be expressed as
\begin{gather}
\begin{bmatrix}
t_0         & (\tkk)^{\ast}    \\
\tkk        &      \Tkk        
\end{bmatrix}
\begin{bmatrix}
\vkm \\
\vkk 
\end{bmatrix}
 =
\begin{bmatrix}
\bkm \\ 
\bkk
\end{bmatrix}
\label{evenstep}
\end{gather}
with $c^{(\ell+2)} = [\vkm, (\vkk)^T ]^T = \cmm$. Observe that \eqref{evenstep}
cannot be transformed to a system of the form \eqref{oddstep} by
simple permutations, i.e.\ just by interchanging rows and columns. 
Since we have already solved the systems $\Tcbkk$ and $\Tkk \ykk = - \tkk$
we can write
\begin{equation}
\vkk =(\Tkk)^{-1} (\bkk - \tkk \vkm)=\ckk + \vkm \ykk 
\notag
\end{equation}
and
\begin{align}
\vkm = & (\bkm - [\tkk]^{\ast} \vkk)/t_0 \notag \\
     = & (\bkm-[\tkk]^{\ast} \ckk -[\tkk]^{\ast} \vkm \ykk)/t_0 \notag \\
     = & (\bkm - [\tkk]^{\ast} \ckk)/\betakkk \,, \notag 
\end{align}
where we have used in the last step that $\Tk = [\Tk]^{\ast}$ which implies
that $(\tkk)^{\ast} \ykk$ is real and therefore 
$t_0+(\tkk)^{\ast} \ykk=t_0 +(\tkk)^{T}\conj{\ykk}=\betakkk$.

Note that at each level $M$ we have to check if the stopping criterion
\eqref{stop} is satisfied. The evaluation of the expression
\begin{equation}
\sum_{j=1}^{r} |\pmlsp (x_j) - \smp_j|^2 w_j
\label{lserror}
\end{equation}
can be considerably simplified and by avoiding the evaluation of $\pmlsp$
at the nonuniformly spaced points $x_j$ we can reduce the computational 
effort from $\ord(Mr)$ to $\ord(M)$ operations. 

To do this we define the subspace 
${\cal R} = \big\{ \{p(x_j)\}_{j=1}^{r} :
p \in \PM\big\} \subseteq \Cst^r$ with the weighted inner product
$\langle y, z \rangle_{{\cal R}} = \sum_{j=1}^{r} y_j \bar{z}_j w_j$ for
$y,z \in \Cst^r$. The solution of the least squares problem \eqref{LSP} is
the orthogonal projection of the vector $\{\smp_j\}_{j=1}^{r} \in \Cst^r$ onto
${\cal R}$ and therefore must satisfy
$$
\langle \{\pmlsp (x_j)\} - \{\smp_j\},\{\pmlsp (x_j)\}\rangle_{\cal R} 
= \sum_{j=1}^{r} (\pmlsp (x_j)- \smp_j)\overline{\pmlsp (x_j)} w_j = 0
$$
which implies
\begin{equation}
\langle \{\pmlsp (x_j)\} , \{\smp_j\}\rangle_{\cal R} =
\langle \{\pmlsp (x_j)\} , \{\pmlsp(x_j)\}\rangle_{\cal R}  \,.
\label{ortho}
\end{equation}
Since
\begin{align}
\sum_{j=1}^{r} |\smp_j - \pmlsp (x_j)|^2 w_j = & \sum_{j=1}^{r}|\smp_j|^2 w_j -
2\,\mbox{Re}\,\langle \smp,\{ \pmlsp (x_j)\}\rangle_{\cal R}+
\sum_{j=1}^{r} | \pmlsp (x_j)|^2 w_j \notag \\
=& \sum_{j=1}^{r}|\smp_j|^2 w_j - \sum_{j=1}^{r} |\pmlsp (x_j)|^2 w_j \notag
\end{align}
by~\eqref{ortho}, and because
\begin{align}
\sum_{j=1}^{r} |\pmlsp (x_j)|^2 w_j &= \sum_{j=1}^{r} w_j 
\bigg( \sum_{m=-M}^{M} \cm_m e^{2 \pi i m x_j} \bigg)
\bigg( \sum_{n=-M}^{M} \overline{\cm_n} e^{2 \pi i n x_j} \bigg) \notag \\
&=\sum_{m=-M}^{M}\sum_{n=-M}^{M} \cm_m \overline{\cm_n} 
\bigg( \sum_{j=1}^{r} w_j e^{2 \pi i(m-n)x_j} \bigg) \notag \\
&= \langle \TM \cm , \cm  \rangle = \langle \bm, \cm \rangle \,,
\end{align}
it follows that
\begin{equation}
\sum_{j=1}^{r} |\smp_j - \pmlsp (x_j)|^2 w_j = \sum_{j=1}^{r}|\smp_j|^2 w_j
-\langle \bm,\cm \rangle \,.
\label{stopsimp}
\end{equation}
Since $\sum_{j=1}^{r}|\smp_j|^2 w_j$ has to be computed only once at the 
beginning of the algorithm, the evaluation of \eqref{lserror} can be 
carried out in $\ord(M)$ operations. 

Summing up we have arrived at the following algorithm to compute $\pnstop$.

\medskip

\begin{algorithm}[Levinson-Galerkin algorithm for trigonometric polynomials]
\label{alg1} 
Let the sampling points $\{x_j\}_{j=1}^{r}$, sampling values 
$\{\smp_j\}_{j=1}^{r}$, weights $w_j>0$ and the data error estimate
$\eps$ be given.
Then the trigonometric polynomial $\pnstop$ determined in Algorithm~\ref{alg}
can be computed in $\ord(r \Nstop + \Nstop^2)$ operations by the
following algorithm.\\ 

\noindent
Initialize: 
$t_0 = \sum_{j=1}^{r} w_j, t_1 = \sum_{j=1}^{r} w_j e^{2\pi i x_j},
b_0 = \sum_{j=1}^{r} \smp_j w_j, \sigma = \sum_{j=1}^{r} |\smp_j|^2 w_j$,
$y^{(1)}=-t_1/t_0, c^{(1)}= b_0/t_0, \beta_0=t_0, \alpha_0=-t_1/t_0, 
\eps_1 = (\sigma - b_0^2/t_0)/\sigma, \ell=1$.

\allowdisplaybreaks
\begin{eqnarray}
  &  & \mbox{\bf while} \ \, \eps_\ell > \eps \nonumber \notag \\ 
  &  & \qquad \betakk = (1-\alphak \conj{\alphak}) \betak
      \notag\\ 
  &  & \qquad \mbox{\bf if} \ \ell \equiv 1 \!\!\!\mod 2\notag  \\
  &  & \qquad \qquad \bkl= \sum_{j=1}^{r} \smp_j w_j e^{\pi i (\ell+1) x_j} 
       \notag\\
  &  & \qquad \qquad \vkl=
       \frac{\bkl-\langle \Ek \overline{c^{(\ell)}}, \tk \rangle }
                                                {\betakk}\notag  \\
  &  & \qquad \qquad \vk= \ck+\vkl \Ek \overline{y^{(\ell)}}
       \notag  \\
  &  & \qquad \qquad \ckk = 
              \begin{bmatrix} \vk \\ \vkl \end{bmatrix}\notag\\
  &  & \qquad \qquad \bkk =
              \begin{bmatrix} \bk \\ \bkl \end{bmatrix}\notag  \\
  &  & \qquad {\bf elseif} \ \, \ell \equiv 0\!\!\! \mod 2\notag  \\  
  &  & \qquad \qquad \bkn= \sum_{j=1}^{r} \smp_j w_j e^{-\pi i \ell x_j} 
       \notag\\
  &  & \qquad \qquad \vkn=
       \frac{\bkn-\langle \overline{c^{(\ell)}}, \tk \rangle}
                                                  {\betakk}\notag  \\
  &  & \qquad \qquad \vk = \ck + \vkn \yk \notag  \\
  &  & \qquad \qquad \ckk = 
              \begin{bmatrix} \vkn \\ \vk \end{bmatrix}\notag  \\
  &  & \qquad \qquad \bkk = 
              \begin{bmatrix} \bkn\\ \bk \end{bmatrix}\notag  \\
  &  & \qquad \qquad \eps_{\ell+1} =|\sigma - \langle \bkk,\ckk \rangle|/\sigma \notag \\
  &  & \qquad \mbox{\bf end}\notag  \\
  &  & \qquad \tkk = \sum_{j=1}^{r} w_j e^{2\pi i (\ell+1) x_j} \notag\\
  &  & \qquad \alphakk = 
  - \frac{\tkk+\langle \Ek \overline{y^{(\ell)}},\tk \rangle }
                                           {\betakk}\notag  \\
  &  & \qquad z^{(\ell)}=\yk+\alphakk \Ek \overline{y^{(\ell)}}
   \notag \\
  &  & \qquad \ykk = 
       \begin{bmatrix} z^{(\ell)} \\ \alphakk \end{bmatrix}\notag \\
  &  & \qquad \tkk =  
            \begin{bmatrix} \tk \\ t_{\ell+1} \end{bmatrix}\notag \\
  &  & \qquad \ell = \ell +1 \notag \\
  &  & \mbox{\bf end} \notag\\
  &  & \Nstop = \ell/2\notag \\
  &  & \pnstop(x) = \sum_{k=-\Nstop}^{\Nstop} \cnstop_k e^{2\pi i kx} \notag
\end{eqnarray}
\end{algorithm}

\remark Usually one evaluates the final approximation on regularly spaced grid
points, hence the last step of the algorithm can be realized by a Fast
Fourier transform. The most costly steps are the computation of the
entries of $\tk$ and $\bk$. According to Corollary~1 in \cite{FGS95} the
entries of $\TM$ and $\bm$ can also be computed via FFT by embedding
the $x_j$ into a regular grid (since the $x_j$ can be stored only in finite
precision). In this case one automatically gets all
entries $t_0, \dots, t_{r}$ at once. However this trick is only useful if
the number of points of the regular grid is of the same magnitude as
the number of sampling points. Alternatively one may use
the numerical attractive formulas of Rokhlin 
\cite{DR93} or Beylkin \cite{Bey95} for a fast evaluation of
trigonometric sums at unequally spaced nodes.

Algorithm~\ref{alg1} can be simplified for real-valued data, 
this modification is left to the reader.

Fast Vandermonde solvers require $\ord(Mr)$ operations for the solution of
$\VM \cm = \smp$, cf.~\cite{RAG91}. It is not clear however if these algorithms can 
utilize the nested structure of the sequence of matrices $\{\VM\}_M$ in order to
give rise to an efficient implementation of Algorithm~\ref{alg}. Moreover it is
an open problem if the Vandermonde solvers can be extended to multivariate
trigonometric approximation. We will see in the next section that the extension
of Algorithm~\ref{alg1} to higher dimensions is straightforward.

\section{Multivariate trigonometric approximation} 
\label{s:higherdim}
 
An advantage of the proposed approach, besides its numerical efficiency,
is the fact that it can be easily extended to multivariate trigonometric
approximation. 
In this section we briefly discuss some results for the 
two-dimensional case.

We define the space of 2-D trigonometric polynomials $\PMM$ by
\begin{equation}
\PMM = \left\{ p : p(x,y) = \sum_{j,k=-M}^{M} c_{j,k} e^{2\pi i (jx+ky)} 
\right\}\,.
\label{pm2}
\end{equation}
To reduce the notational burden, we have assumed in~\eqref{pm2} that 
$p$ has equal degree $M$ in each coordinate, the extension to polynomials 
with different degree in each coordinate is straightforward.

For an arbitrary sampling set $\{(x_j,y_j)\}_{j=1}^{r} 
\in [0,1)^2$ and given degree $M$ the system matrix according to
the 2-D version of Theorem~\ref{thm1} is~\cite{Str97} 
\begin{equation}
(\TM)_{k,l} = \sum_{j=1}^{r} w_j e^{2 \pi i (k-l)(x_j+y_j)}\,, \quad
k,l = 0,\dots 2M\,.
\label{bttb}
\end{equation}
One can easily verify that $\TM$ is a hermitian block Toeplitz matrix with 
$2M+1$ different Toeplitz blocks of size $(2M+1) \times 2M+1$, 
cf.~\cite{Str97}. 
For a given sampling set let $\TM$ be the block Toeplitz matrix for
degree $M$ and $\TMM$ the block Toeplitz matrix for degree $M+1$.
There is a similar relationship between $\TM$ and $\TMM$ as in
the 1-D case. More precisely, denote the Toeplitz blocks of $\TM$ and $\TMM$
by $(\BKM)_{k}, k=0,\dots,2M$, and $(\BKMM)_{k}, k=0,\dots,2M+2$, respectively.
Then one readily verifies the following embedding:
\allowdisplaybreaks
\begin{gather}
\label{blockembed}
\TMM = 
\begin{bmatrix}
(\BKMM)_0 & \dots & (\BKMM)^{\ast}_{2(M+1)} \\
\vdots & \fbox{$\TM$} & \vdots \\
(\BKMM)_{2(M+1)} & \dots & (\BKMM)_0
\end{bmatrix},\\
\BKMM = 
\begin{bmatrix}
t_0 & \dots & \conj{t_{2(M+1)}} \\
\vdots & \fbox{$\BKM$} & \vdots \\
t_{2(M+1)} & \dots & t_0
\end{bmatrix}.
\end{gather}

In \cite{Aka73} Levinson's algorithm has been extended to general
block Toeplitz systems. With this extension and relation~\eqref{blockembed} 
at hand, we can easily generalize Algorithm~\ref{alg1} to 2-D
(and along the same lines to multivariate) trigonometric approximation. 

The analysis of the stopping criterion~\eqref{stop} in Section~\ref{s:multi}
can be applied line by line to the 2-D (actually to the $n$-D) setting.
The only difficulty arises in the search for simple criteria for the
invertibility of the block Toeplitz matrix $\TM$. The condition 
$$(2M+1)^d \le r$$
is necessary in dimension $d>1$, but no longer sufficient, since the 
fundamental theorem of algebra does not hold in dimensions larger than one. 
In \cite{Gro97} Gr\"ochenig has derived estimates for the condition number 
of $\TM$ in higher dimensions. In 2-D these estimates can be stated as follows.

Let $D_{\delta } (a,b) = \{ (x,y) \in \Rst ^2: (x-a)^2+(y-b)^2 < \delta ^2 \}$ 
be the disc of radius $\delta $ centered at $(a,b)$. We say that a set
$\{(x_j,y_j), j=1, \dots ,r\} $ is $\delta$-dense in $[0,1]\times [0,1]$, if
$\bigcup _{j=1}^r D_{\delta }(x_j,y_j) \supseteq  [0,1]\times [0,1] $.
In other words, the distance of a given sample $(x_j,y_j)$ to its
nearest neighbor $(x_k,y_k), k\neq j$ is at most $2 \delta $.

Analogously to Section~\ref{s:pred} we choose the size of the Voronoi region
$V_j$ associated with $x_j$ as weight $w_j$ in the computation of the
block Toeplitz matrix $\TM$ in~\eqref{bttb}.
Suppose that the sampling set $\{(x_j, y_j), j=1 , \dots ,r \}
\subseteq [0,1]\times [0,1]$ is $\delta $-dense and
\begin{equation}
\delta < \frac{\log 2}{4\pi M } \, .
\label{density}
\end{equation}
Gr\"ochenig~\cite{Gro97} has shown that under these conditions
\begin{equation}
\label{blockcond}
\cond(\TM) \le \frac{4}{(2 - e^{4 \pi M\delta })^2}.
\end{equation}
In particular, for arbitrary $\delta $-dense sampling sets, the block
Toeplitz matrix $\TM$ is invertible and the 2-D version of
Algorithm~\ref{alg1} is applicable.

\subsection{Line-type nonuniform sampling in 2-D} \label{ss:line}

In the following we consider a special case of trigonometric approximation
in two dimensions. This case arises when a function is irregularly 
sampled along lines. A typical example is illustrated in 
Figure~\ref{fig:line}. Such sampling patterns are encountered for instance 
in geophysics and medical imaging, see also Section~\ref{ss:seq}.

\begin{figure}[h]
\begin{center}
\epsfig{file=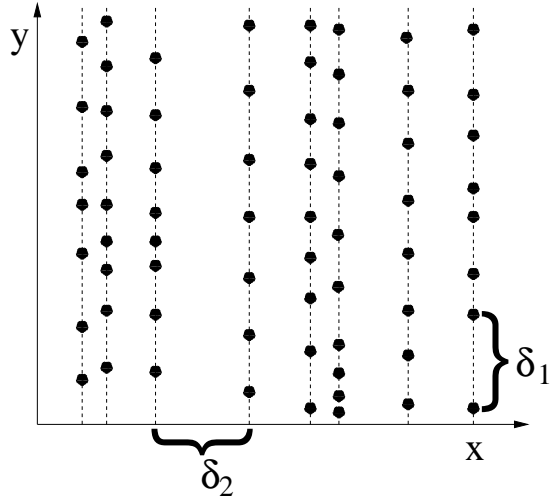,width=70mm}
\caption{Line-type nonuniform sampling set}
\label{fig:line}
\end{center}
\end{figure}

\begin{corollary} \label{cor1}
Let $p \in \PMM$ and let $\{x_j,y_{j,k}\}, j=1,\dots,r,
k=1,\dots,r_j$ be a sampling set in $[0,1)^2$ such that
\begin{equation}
A_1 \|p\|^2 \le \sum_{k=1}^{r_j} |p(y_{j,k})|^2 \le B_1 \|p\|^2
\qquad A_1, B_1 >0
\label{ystable}
\end{equation}
for every $p \in \PM$ and for all $j$.
Further assume that $\{x_j\}_{j \in \Zst}$ is a sampling set such that
\begin{equation}
A_2 \|p\|^2 \le \sum_{j=1}^{r} |p(x_{j})|^2 \le B_2 \|p\|^2
\qquad A_2, B_2 >0
\label{xstable}
\end{equation}
for every $p \in \PM$. Then
\begin{equation}
A_1 A_2 \|p\|^2 \le \sum_{j=1}^{r} \sum_{k =1}^{r_j}
|p(x_j,y_{j,k})|^2 \le B_1 B_2 \|p\|^2
\label{xystable}
\end{equation}
for every $p \in \PMM$.

If $\{x_j\}$ and $\{y_{j,,k}\}$ are sampling sets with
$\sup_{k} (x_{k+1} - x_{k}) = \delta_2 < \frac{1}{2M}$ and
$\sup_{j,k} (y_{j,k+1}-y_{j,k}) = \delta_1 < \frac{1}{2M}$,
then $A_l = (1-\delta_l)^2, B_l = (1+\delta_l)^2, l=1,2$ and 
the condition number of the block Toeplitz matrix $\TM$ is bounded by
\begin{equation}
\kappa(\TM) \le \frac{(1+\delta_1)^2 (1+\delta_2)^2}
                   {(1-\delta_1)^2 (1-\delta_2)^2} \,.
\label{btoepcond}
\end{equation}
\end{corollary}
\begin{proof} 
Let $x$ be fixed. Then $y \rightarrow p(x,y) \in \PM$ and hence
for all $j$
\begin{equation}
A_1 \int \limits_{0}^{1} |p(x_j,y)|^2 \, dy \le
                      \sum_{k=1}^{r_j} |p(x_j,y_{j,k})|^2 \le 
B_1 \int \limits_{0}^{1} |p(x_j,y)|^2 \, dy 
\label{ysum}
\end{equation}
by assumption~\eqref{ystable}.
It follows that
\begin{equation}
A_1 \sum_{j=1}^{r} \int \limits_{0}^{1} |p(x_j,y)|^2 \, dy \le
\sum_{j,k} |p(x_j,y_{j,k})|^2 \le 
B_1 \sum_{j=1}^{r} \int \limits_{0}^{1}  |p(x_j,y)|^2 \, dy 
\label{xysum}
\end{equation}
Now let $y$ be fixed. Then $x \rightarrow p(x,y) \in \PM$ and 
\begin{equation}
A_2 \int \limits_{0}^{1} |p(x,y)|^2 \, dx 
\le \sum_{j=1}^{r} |p(x_j,y)|^2 \le
B_2 \int \limits_{0}^{1} |p(x,y)|^2 \, dx  \,.
\label{xsum}
\end{equation}
Since 
\begin{equation}
\sum_{j=1}^{r} \int \limits_{0}^{1}  |p(x_j,y)|^2 \, dy 
 = \int \limits_{0}^{1} \sum_{j=1}^{r} |p(x_j,y)|^2 \, dy\,, 
\label{fubini}
\end{equation}
assertion \eqref{xystable} follows by combining \eqref{xysum} 
and \eqref{xsum} with~\eqref{fubini}. The estimate of the constants 
$A_l,B_l$ and of the condition number of the block Toeplitz matrix $\TM$ 
follow from Theorem~\ref{thm1}.
\end{proof}

The proof of Corollary~\ref{cor1} is due to Gr\"ochenig~\cite{Gro}.
Corollary~\ref{cor1} does not only guarantee that $p \in \PMM$ can be 
recovered from its samples $p(x_j,y_{j,k})$, it 
provides more. An immediate consequence is, that it can be reconstructed
by an efficient algorithm relying on an successive application of 
Algorithm~\ref{alg1} and the Gohberg-Semencul representation
of the inverse of a Toeplitz matrix. See Section~\ref{ss:seq} for more
details and an application in medical imaging.

\section{Curve and surface approximation by trigonometric polynomials}
\label{s:curve}

Trigonometric polynomials can be used to model the boundary 
or the surface of smooth objects.
Let us consider a two-dimensional object, obtained e.g.\ by a planar
cross-section from a 3-D object and assume that the boundary of this
2-D object is a closed curve in $\Rst^2$. We denote this curve by $f$ and 
parameterize it by $f(u) = (x_u,y_u)$, where $x_u$ 
and $y_u$ are the coordinates of $f$ at ``time'' $u$ in the $x$- and 
$y$-direction respectively. Obviously we can interpret $f$ as a 
one-dimensional continuous, complex, and periodic
function, where $x_u$ represents the real part and $y_u$ represents the
imaginary part of $f(u)$. It follows from the Theorem of Weierstrass 
(and from the Theorem of Stone-Weierstrass~\cite{Rud76} for higher dimensions) 
that a continuous
periodic function can be approximated uniformly by trigonometric polynomials.
If $f$ is smooth, we can fairly assume that trigonometric polynomials of
low degree provide an approximation of sufficient precision.

Assume that we know only some arbitrary, perturbed points 
$s_j= (x_{u_j},y_{u_j})=f(u_j)+\delta_j, j=1,\dots, r$ of $f$, and we want to 
recover $f$ from these points. By a slight abuse of notation we interpret
$s_j$ as complex number and write
\begin{equation}
s_j = x_j + i y_j\,.
\label{transsamples}
\end{equation}
We relate the curve parameter $u$ to the boundary points $s_j$ by 
computing the distance between two successive points $s_{j-1},s_j$ via
\begin{align}
\label{transpoints1}
u_1 & = 0 \\
\label{transpoints2}
u_j & = u_{j-1} + d_j \\
\label{transpoints3}
d_j & = \sqrt{(x_j - x_{j-1})^2 + (y_j - y_{j-1})^2}
\end{align}
for $j=2,\ldots,r$. Via the normalization $t_j = t_j/L$ with
$L = u_{r} + d_N$ we force all sampling points to be in $[0,1)$.
Other choices for $d_j$ in~\eqref{transsamples} can be found 
in~\cite{Die93} in conjunction with curve approximation using splines. 

Having carried out the transformations 
\eqref{transsamples}--\eqref{transpoints3}, we can solve the problem of 
recovering the curve $f$ from its perturbed points $s_j$ by 
Algorithm~\ref{alg1}.

\subsection{Object boundary recovery in Echocardiography} \label{ss:echo}

Trigonometric polynomials are certainly not suitable to model
the shape of arbitrary objects. However they are often useful in cases
where an underlying (stationary) physical process implies smoothness 
conditions of the object. Typical examples arise in medical imaging, for 
instance in clinical cardiac studies, where the evaluation of cardiac
function using parameters of left ventricular contractibility is an important 
constituent of an echocardiographic examination~\cite{WGL93}. These 
parameters are derived using boundary tracing of endocardial borders
of the Left Ventricle (LV). The extraction of the boundary of the LV 
comprises two steps, once the ultrasound image of a cross section of the
LV is given, see Figure~\ref{fig:cardio}(a)--(d). First an edge detection 
is applied to the ultrasound image to detect the boundary of the
LV, cf.~Figure~\ref{fig:cardio}(c). However this procedure
may be hampered 
by the presence of interfering biological structures (such as papillar 
muscles), the unevenness of boundary contrast, and various kinds of
noise \cite{SBS95}. Thus edge detection often provides only a set of 
nonuniformly spaced, perturbed boundary points rather than a connected boundary. 
Therefore a second step is required, to recover the original boundary from 
the detected edge points, cf.\ Figure~\ref{fig:cardio}(d). Since the shape
of the Left Ventricle is definitely smooth, trigonometric polynomials 
are particularly well suited to model its boundary.  
\begin{figure}
\begin{center}
\subfigure[2-D echocardiography]{
\epsfig{file=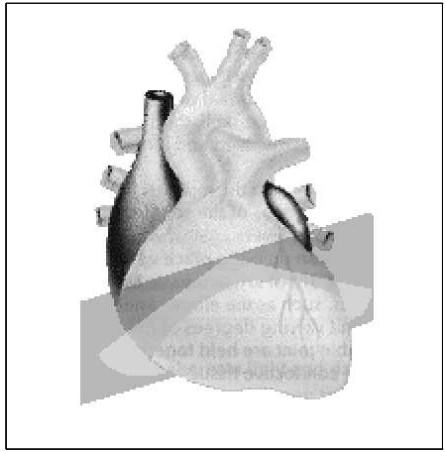,width=60mm}} \quad
\subfigure[Cross section of Left Ventricle]{
\epsfig{file=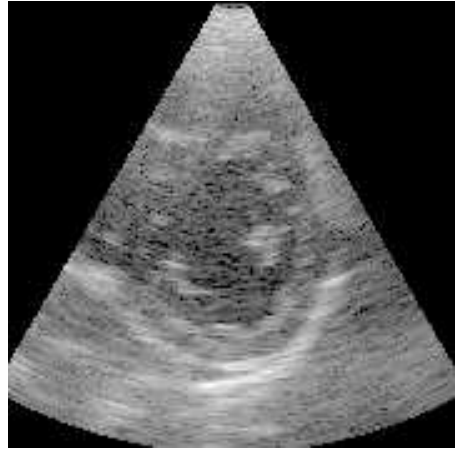,width=60mm}} \\
\vspace*{2mm}
\subfigure[Detected boundary points]{
\epsfig{file=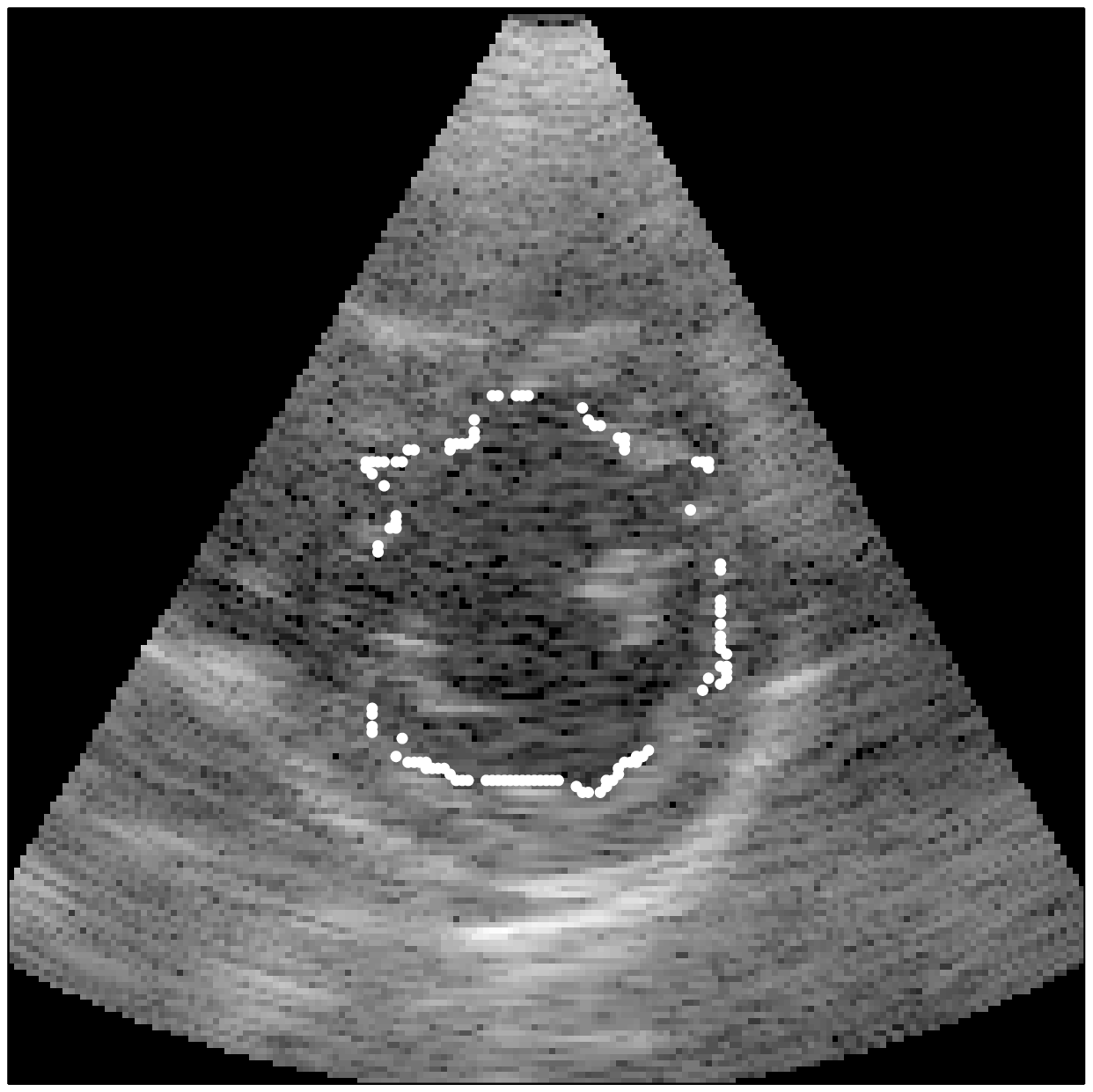,width=60mm}} \quad
\subfigure[Recovered boundary of LV computed by Algorithm~\ref{alg1}]{
\epsfig{file=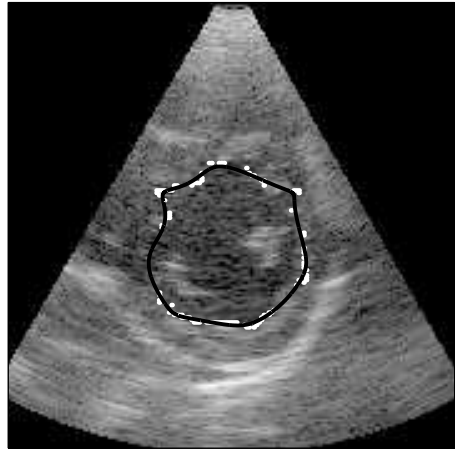,width=60mm}}
\caption{The recovery of the boundary of the Left Ventricle from 
2-D ultrasound images is a basic step 
in echocardiography to extract relevant parameters of cardiac function.}
\label{fig:cardio}
\end{center}
\end{figure}

After having transformed the detected boundary points as described 
in \eqref{transsamples}--\eqref{transpoints3} we can
use Algorithm~\ref{alg1} to recover the boundary. The noise
level $\delta$ depends on the technical equipment under use, it can
be determined from experimental experience.
Figure~\ref{fig:cardiorec}(a)--(b) demonstrate the importance of
determining a proper degree for the approximating polynomial. The 
approximation displayed in Figure~\ref{fig:cardiorec}(a) has been
computed by solving \eqref{LSP} where $M$ has been chosen too small,
we obviously have underfitted the data. The overfitted approximation 
obtained by solving \eqref{LSP} using a too large $M$ is shown in 
Figure~\ref{fig:cardiorec}(b). The
approximation shown in Figure~\ref{fig:cardio}(d) has been computed
by Algorithm~\ref{alg1}, it provides the optimal balance between 
fitting the data and smoothness of the solution. 

\begin{figure}
\begin{center}
\subfigure[Underfitted solution]{
\epsfig{file=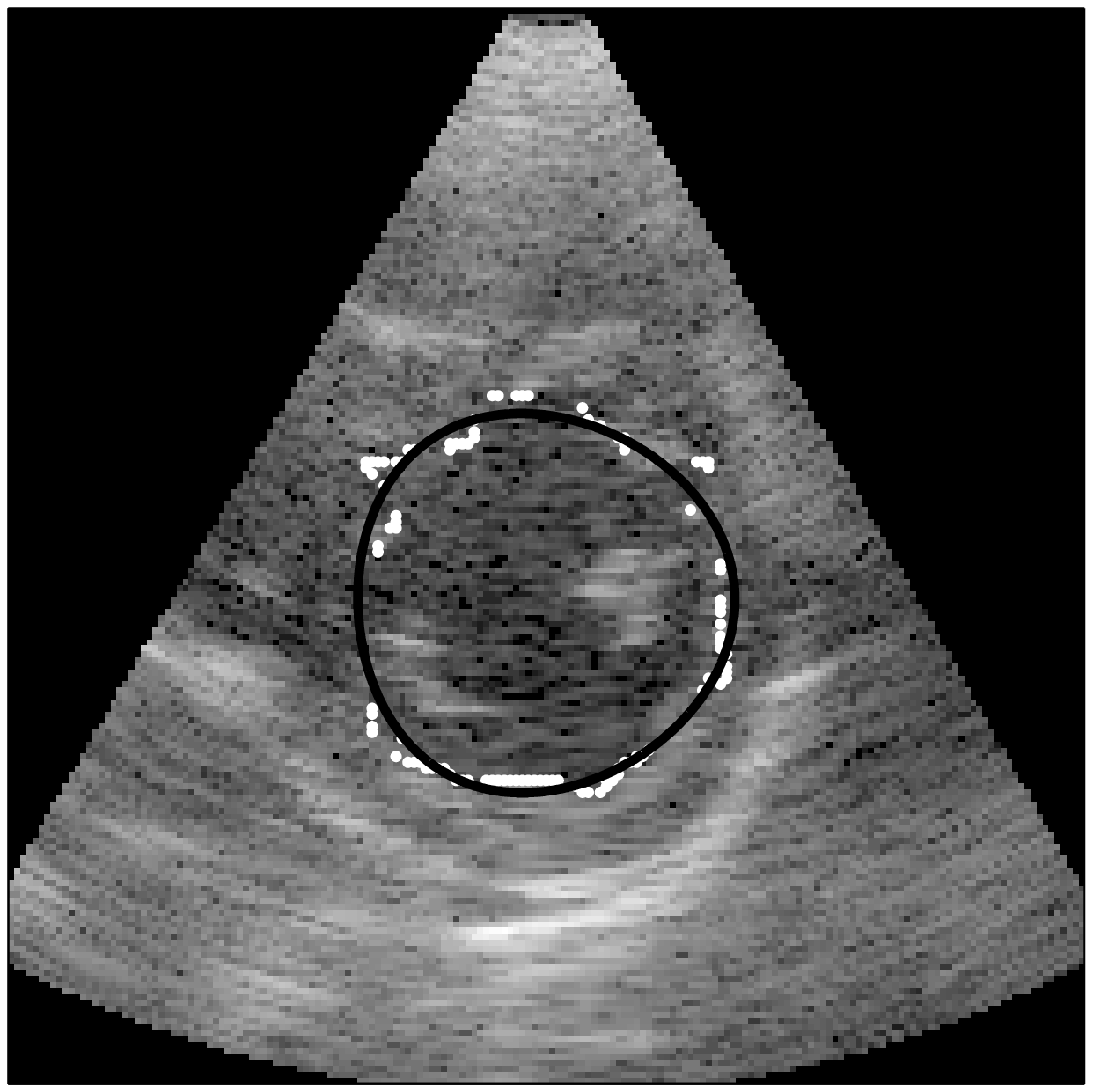,width=60mm}} \quad
\subfigure[Overfitted solution]{
\epsfig{file=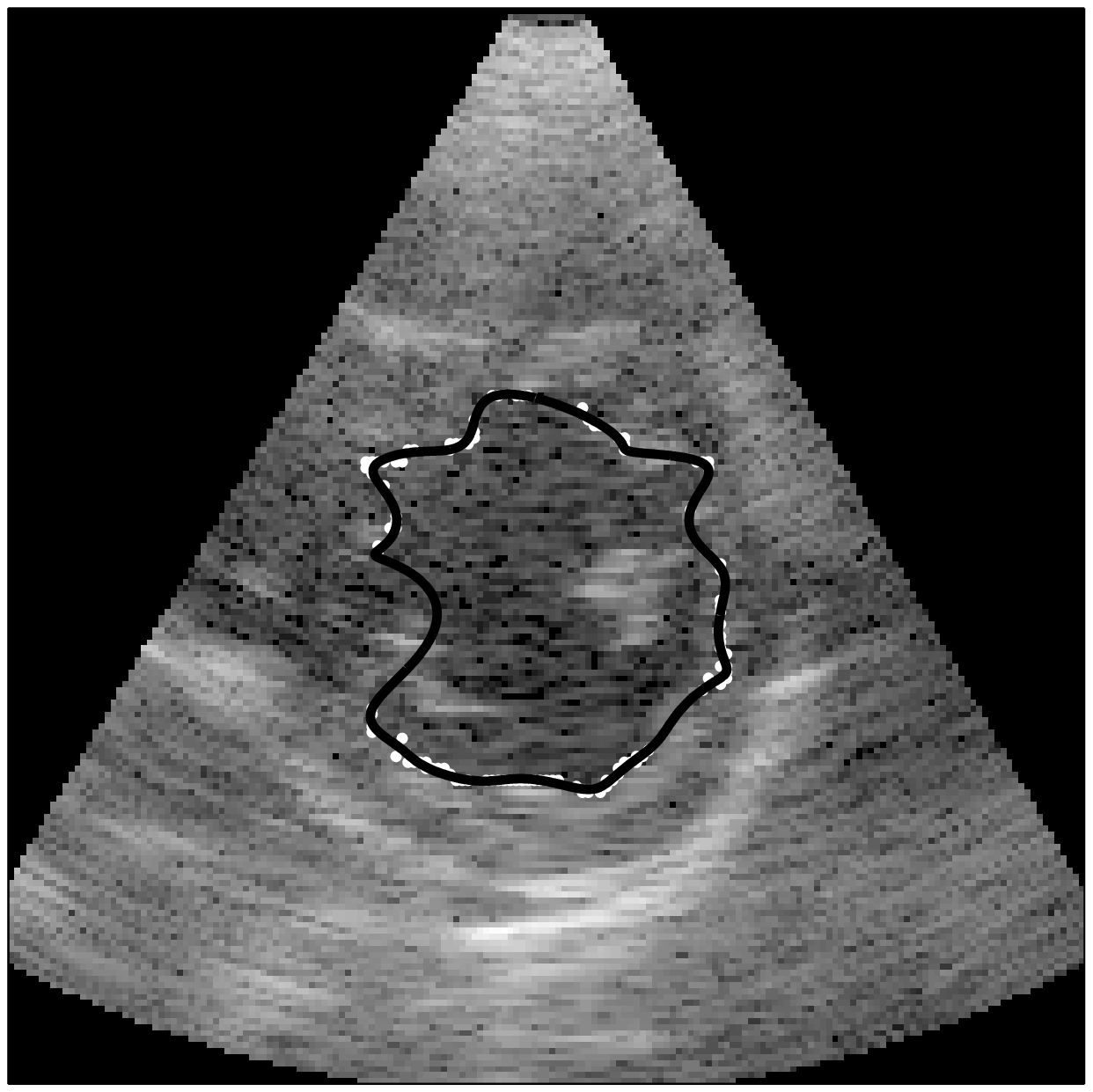,width=60mm}} \\
\caption{The approximation in the left image results from using a 
too small polynomial degree, the approximation in the right image from
a too large degree for the trigonometric approximation.} 
\label{fig:cardiorec}
\end{center}
\end{figure}

\subsection{Boundary recovery from a sequence of images} \label{ss:seq}

In cardiac clinical studies one is more interested in 
the behavior of the Left Ventricle over a period of time rather
than in a single ``snapshot''. Thus for a fixed cross section
we are given a sequence of ultrasound images (usually regularly spaced 
in time) describing the variation of the shape of the LV with time.
One cycle from diastole (the state of maximal contraction of the LV),
passing systole (the state of maximal expansion) to the next diastole 
consists typically of about 30 image frames. Since the behavior of the LV 
is (at least for a short period of time) almost periodic, one can model
the varying shape of a fixed cross section of the LV as distorted 
two-dimensional torus, which in turn can be interpreted as 2-D
trigonometric polynomial. Clearly we have to use a different degree for
the time coordinate $\tau$ and for the spatial coordinate $u$. 

Due to interfering biological structures and other distortions
it sometimes happens that some of the image frames cannot be used to
extract any reliable boundary information. Thus we have to approximate these
missing boundaries from the information of the other image frames.
To be more precise, assume that an echocardiographic examination provides a 
sequence of ultrasound images $I_\tau$ taken at time points 
$\tau=0,1,\dots,T-1$, where $T$ is approximately the length of one
diastolic cycle (the time points could also be nonuniformly spaced).
Assume that some of the images $I_\tau$ provide no useful information, 
so that we can only detect boundary 
points $\{s_{j,k}\}_{k=1}^{r_j}$ from the images $I_{\tau_j}$, where
$\{\tau_j\}_{j=1}^{r}$ is a subset of $0,1,\dots T-1$. In order to
get a complete description of the LV for the time interval $[0,T]$, 
we have not only to
approximate the boundaries $f_j$ from each $I_j$, but we also
have to recover the boundaries corresponding to the missing images.
In other words we look for a 2-D trigonometric polynomial $\popt \in \PMM$ of
appropriate degree $M$ that satisfies 
$p(\tau_j,u_{j,k}) \equiv (x_{j,k},y_{j,k})$
where the parameter $u$ is related to $s_{j,k} = x_{j,k} + i y_{j,k}$ by
formulas~\eqref{transpoints1}--\eqref{transpoints3}. This approximation 
can be computed by the 2-D version of Algorithm~\ref{alg1}, as indicated 
in the beginning of Section~\ref{s:higherdim}.

Under certain conditions we can use the 1-D version of
Algorithm~\ref{alg1} instead of its 2-D version.
As long as the assumptions of Corollary~\ref{cor1} are satisfied,
we can compute $\popt \in \PMM$ by a successive application of
Algorithm~\ref{alg1}. We first approximate the boundaries $f_j$ for each
$j$ separately from its samples $\{s_{j,k}\}_{k=1}^{r_j}$, which 
yields $j$ different polynomials $p^{(M_j)} \in \PMj$. 
Having done this, the next step is to recover the missing boundaries at those
time points where no information is available. 
We proceed by approximating successively the missing information
``line by line''. We choose $u=0$, say, and approximate
the missing information from the samples $p^{(M_j)}(u)$ taken at the time
points $\tau_j, j=1,\dots,r$. 

Note that the Toeplitz matrices of 
the systems $(\TM)_u \cm_u = \bm_u$ coincide for all $u$, since the 
sampling geometry is constant along the $u$-coordinate (because we have 
recovered all samples at each $\tau_j$). Thus we have to solve multiple 
Toeplitz systems with the same system matrix but different right hand side. 
It is well-known that
this can be done efficiently by exploiting the Gohberg-Semencul
representation of the inverse of the Toeplitz matrix~\cite{GS72}. 
In our context this reads as follows. We solve 
\begin{equation}
(\TM)_u \cm_u  = \bm_u
\label{gs1}
\end{equation}
for one $u$ by Algorithm~\ref{alg1}. We can solve now all other systems 
efficiently by establishing $(\TM)^{-1}$ in the Gohberg-Semencul form
\begin{equation}
(\TM)^{-1} = \left([L^{(M)}]^\ast L^{(M)} - U^{(M)} [U^{(M)}]^\ast\right)/z_0 
\label{gs2}
\end{equation}
where $L^{(M)}$ is a lower triangular Toeplitz matrix with 
$z = [z_0, z_1,\dots,z_{2M}]^T$ as its first column, 
$U^{(M)}$ is an upper triangular Toeplitz matrix with 
$[0,z_{1}, \dots, z_{2M}]^T$ as its last column, $z$ being the first
column of $(\TM)^{-1}$. The matrix vector multiplications to compute
$\cm_u = (\TM)_u^{-1} \bm_u$ can now be carried out quickly 
using the Fast Fourier transform by embedding $L^{(M)}$ and $U^{(M)}$ into 
circulant matrices.

\section{Miscellaneous remarks} \label{s:con}

For sampling sets with large gaps it can happen that
the system $\Tcb$ gets ill-conditioned with increasing degree $M$
and therefore Algorithm~\ref{alg1} may become unstable \cite{Cyb80}. 
In this case one can use a 
different, more robust approach, which however comes at 
higher computational costs~\cite{SS97}. We solve the system $\TM \cm = \bm$
iteratively, e.g.\ by the conjugate gradient method until a certain
stopping criterion is satisfied at iteration $k$, say, yielding
the solution $\cm_k$. We use this solution as initial guess at the 
next level $M+1$ by setting $\cmm_0 = [0 \,\, (\cm_k)^T \,\, 0]^T$.  
The crucial point in this procedure is to find a
stopping criterion that guarantees convergence of the iterates,
see~\cite{SS97,RS98} for more details.

The computation of the entries of the Toeplitz matrix in 
Section~\ref{s:curve} involves the nodes $u_j$ which in this particular
case depend on the (perturbed) samples $s_j$. Therefore
not only the right hand side $\bm$, but also $\TM$ is subject
to perturbations. Hence in principle one might use the
concept of total least squares (see~\cite{FGH97}) instead of a
least squares approach.
A detailed discussion of this modification is beyond the scope of this paper.


\section*{Acknowledgments}
The major part of this work was completed during my stay at
Stanford University. I want to thank Prof.~David Donoho and the Department 
of Statistics for their hospitality.


\begin{thebibliography}{10}

\bibitem{Aka73}
{\sc H.~Akaike}, {\em Block {T}oeplitz matrix inversion}, SIAM J. Appl. Math.,
  24 (1973), pp.~234--241.

\bibitem{AG88}
{\sc G.~Ammar and W.~Gragg}, {\em Superfast solution of real positive definite
  {T}oeplitz systems}, SIAM J. Matrix Anal. Appl., 9 (1988), pp.~61--76.

\bibitem{Bey95}
{\sc G.~Beylkin}, {\em On the fast {F}ourier transform of functions with
  singularities}, Appl.\ Comp.\ Harm.\ Anal., 2 (1995), pp.~363--381.

\bibitem{BEF94}
{\sc S.~Boyd, L.~El~Ghaoui, E.~Feron, and V.~Balakrishnan}, {\em Linear matrix
  inequalities in system and control theory}, SIAM, Philadelphia, PA, 1994.

\bibitem{CN96}
{\sc R.~Chan and M.~Ng}, {\em Conjugate gradient methods for {T}oeplitz
  systems}, SIAM Review, 38 (1996), pp.~427--482.

\bibitem{Cyb80}
{\sc G.~Cybenko}, {\em The numerical stability of the {L}evinson-{D}urbin
  algorithm for {T}oeplitz systems of equations}, SIAM J.\ Sci.\ Statist.\
  Comp., 1 (1980), pp.~303--3190.

\bibitem{Dem89}
{\sc C.~Demeure}, {\em Fast {QR} factorization of {V}andermonde matrices},
  Linear Algebra Appl., 122--124 (1989), pp.~165--194.

\bibitem{Die93}
{\sc P.~Dierckx}, {\em Curve and surface fitting with splines}, Monographs on
  Numerical Analysis, Oxford University Press, 1993.

\bibitem{DR93}
{\sc A.~Dutt and V.~Rokhlin}, {\em Fast {F}ourier transforms for nonequispaced
  data}, SIAM J.\ Sci.\ Comp., 14 (1993), pp.~1368--1394.

\bibitem{Fas97}
{\sc H.~Fassbender}, {\em On numerical methods for discrete least-squares
  approximation by trigonometric polynomials}, Math.\ Comp., 66 (1997),
  pp.~719--741.

\bibitem{FGS95}
{\sc H.~G. Feichtinger, K.~Gr{\"o}chenig, and T.~Strohmer}, {\em Efficient
  numerical methods in non-uniform sampling theory}, Numerische Mathematik, 69
  (1995), pp.~423--440.

\bibitem{FGH97}
{\sc R.~Fierro, G.~Golub, P.~Hansen, and D.~O'Leary}, {\em Regularization by
  truncated total least squares}, SIAM J.\ Sci.\ Comp., 18 (1997),
  pp.~1223--1241.

\bibitem{Gau91}
{\sc W.~Gautschi}, {\em How (un)stable are {V}andermonde systems?}, in
  Asymptotic and computational analysis (Winnipeg, MB, 1989), Dekker, New York,
  1990, pp.~193--210.

\bibitem{GS72}
{\sc I.~Gohberg and A.~Semencul}, {\em On the inversion of finite {T}oeplitz
  matrices and their continuous analogs}, Mat. Issled., 2 (1972), pp.~201--233.

\bibitem{GHW79}
{\sc G.~Golub, M.~Heath, and G.~Wahba}, {\em Generalized cross-validation as a
  method for choosing a good ridge parameter}, Technometrics, 21 (1979),
  pp.~215--223.

\bibitem{GL96}
{\sc G.~Golub and C.~van Loan}, {\em Matrix Computations}, Johns Hopkins,
  Baltimore, third~ed., 1996.

\bibitem{Gro}
{\sc K.~Gr{\"o}chenig}, {\em personal communication}.

\bibitem{Gro93a}
\leavevmode\vrule height 2pt depth -1.6pt width 23pt, {\em A discrete theory of
  irregular sampling}, Linear Algebra Appl., 193 (1993), pp.~129--150.

\bibitem{Gro97}
\leavevmode\vrule height 2pt depth -1.6pt width 23pt, {\em Finite and
  infinite-dimensional models for non-uniform sampling}, in SampTA - Sampling
  Theory and Applications, Aveiro, Portugal, 1997, pp.~285--290.

\bibitem{Han92}
{\sc P.~Hansen}, {\em Analysis of discrete ill-posed problems by means of the
  {L}-curve}, SIAM Review, 34 (1992), pp.~561--680.

\bibitem{Lev47}
{\sc N.~Levinson}, {\em The {W}iener {RMS} (root-mean square) error criterion
  in filter design and prediction}, J. Math. Phys., 25 (1947), pp.~261--278.

\bibitem{New70}
{\sc A.~Newbery}, {\em Trigonometric interpolation and curve-fitting}, Math.\
  Comp., 24 (1970), pp.~869--876.

\bibitem{OBS92}
{\sc A.~Okabe, B.~Boots, and K.~Sugihara}, {\em Spatial tessellations: concepts
  and applications of {V}orono\u\i\ diagrams}, John Wiley \& Sons Ltd.,
  Chichester, 1992.
\newblock With a foreword by D. G. Kendall.

\bibitem{RS98}
{\sc M.~Rauth and T.~Strohmer}, {\em Smooth approximation of potential fields
  from noisy scattered data}, Geophysics, 63 (1998), pp.~85--94.

\bibitem{RAG91}
{\sc L.~Reichel, G.~Ammar, and W.~Gragg}, {\em Discrete least squares
  approximation by trigonometric polynomials}, Math. Comp., 57 (1991),
  pp.~273--289.

\bibitem{Rud76}
{\sc W.~Rudin}, {\em Fourier Analysis on Groups}, Wiley Interscience, New York,
  1976.

\bibitem{SS97}
{\sc O.~Scherzer and T.~Strohmer}, {\em A multi--level algorithm for the
  solution of moment problems}, Num.Funct.Anal.Opt., 19 (1998), pp.~353--375.

\bibitem{Str97}
{\sc T.~Strohmer}, {\em Computationally attractive reconstruction of
  band-limited images from irregular samples}, IEEE Trans.\ Image Proc., 6
  (1997), pp.~540--548.

\bibitem{SBS95}
{\sc M.~Suessner, M.~Budil, T.~Strohmer, M.~Greher, G.~Porenta, and T.~Binder},
  {\em Contour detection using artifical neuronal network presegmention}, in
  Proc. Computers in Cardiology, Vienna, 1995.

\bibitem{WGL93}
{\sc D.~Wilson, E.~Geiser, and J.~Li}, {\em Feature extraction in 2-dimensional
  short-axis echocardiographic images}, J. Math. Imag. Vision, 3 (1993),
  pp.~285--298.

\end{thebibliography}

\end{document}